\documentclass[11pt,reqno]{amsart}

\usepackage{amsmath}
\usepackage{amsxtra}
\usepackage{amscd}
\usepackage{amsthm}
\usepackage{amsfonts}
\usepackage{amssymb}
\usepackage{eucal}
\usepackage[all]{xy}

\numberwithin{equation}{section}

\newtheorem{thm}{Theorem}[section]
\newtheorem{prop}[thm]{Proposition}
\newtheorem{lem}[thm]{Lemma}
\newtheorem{cor}[thm]{Corollary}
\theoremstyle{definition}
\newtheorem{rem}[thm]{Remark}
\newtheorem{exm}[thm]{Example}
\newtheorem{dfn}[thm]{Definition}
\newcommand{\nc}{\newcommand}
\nc{\C}{{\mathbb C}}

\newcommand{\bean}{\begin{eqnarray*}}
\newcommand{\eean}{\end{eqnarray*}}

\newcommand{\bea}{\begin{eqnarray}}
\newcommand{\ena}{\end{eqnarray}}

\def\ov[#1,#2]{\overset{\scriptstyle #1}{#2}}
\nc{\verm}{M_{k,l}}

\nc{\slt}{{\widehat{\mathfrak{sl}}_2}}
\nc{\inte}{{L_{k,l}}}
\nc{\pkln}{{\mathcal P}_{k,l}^{(N)}}
\nc{\ckln}{{c_{k,l}^{(N)}}}
\def\qbin[#1;#2]{{\left[\matrix{\displaystyle #1}\\{\displaystyle #2}\endmatrix\right]}}
\def\geq{\ge}
\def\leq{\le}

\nc{\Vk}{{\mathfrak V}_k}

\def\xx[#1,#2]{{#1}^{(#2)}}
\def\yy[#1,#2,#3]{{#1}^{(#2)}_{#3}}

\def\v(#1;#2){\Bigl({#1\atop#2}\Bigr)}

\nc{\alb}{{\boldsymbol{\alpha}}}
\nc{\beb}{{\boldsymbol{\beta}}}

\def\hookdownarrow%
{\Big\downarrow\kern-.14cm\smash{\raise.4cm\hbox{$\scriptstyle\cap$}}}

\def\uC{\underline C}
\def\S#1{S[#1,{\bf a}]=l}
\def\L#1{L[#1,{\bf a}]=k+l}
\title[Particle Content]
{Particle content of the $(k,3)$-configurations.}

\author{B. Feigin, M. Jimbo, T. Miwa, E. Mukhin and Y. Takeyama}
\address{BF: Landau institute for Theoretical Physics, Chernogolovka,
142432, Russia}\email{feigin@feigin.mccme.ru}  
\address{MJ: Graduate School of Mathematical Sciences, University of
Tokyo, Tokyo 153-8914, Japan}\email{jimbomic@ms.u-tokyo.ac.jp}
\address{TM: Department of Mathematics, Graduate School of Science, 
Kyoto University, Kyoto 606-8502 Japan}\email{tetsuji@kusm.kyoto-u.ac.jp}
\address{EM: Department of Mathematics, 
Indiana University-Purdue University-Indianapolis, 
402 N.Blackford St., LD 270, 
Indianapolis, IN 46202}
\email{mukhin@math.iupui.edu}
\address{YT: Department of Mathematics, Graduate School of Science, 
Kyoto University, Kyoto 606-8502 Japan}
\email{takeyama@kusm.kyoto-u.ac.jp} 

\date{\today}
\begin{document}
\setcounter{footnote}{0}\renewcommand{\thefootnote}{\arabic{footnote}}

\begin{abstract}
For all $k$, we construct a bijection between the set of sequences
of non-negative integers ${\bf a}=(a_i)_{i\in{\bf Z}_{\geq0}}$ satisfying
$a_i+a_{i+1}+a_{i+2}\leq k$ and the set of rigged partitions $(\lambda,\rho)$.
Here $\lambda=(\lambda_1,\ldots,\lambda_n)$ is a partition
satisfying $k\geq\lambda_1\geq\cdots\geq\lambda_n\geq1$
and $\rho=(\rho_1,\ldots,\rho_n)\in{\bf Z}_{\geq0}^n$ is such that
$\rho_j\geq\rho_{j+1}$ if $\lambda_j=\lambda_{j+1}$. One can think of $\lambda$
as the particle content of the configuration ${\bf a}$ and $\rho_j$ as
the energy level of the $j$-th particle, which has the weight $\lambda_j$.
The total energy $\sum_iia_i$ is written as the sum of the
two-body interaction term $\sum_{j<j'}A_{\lambda_j,\lambda_{j'}}$
and the free part $\sum_j\rho_j$. The bijection implies
a fermionic formula for the one-dimensional configuration sums
$\sum_{\bf a}q^{\sum_iia_i}$. We also derive the polynomial identities which
describe the configuration sums corresponding to the configurations with
prescribed values for $a_0$ and $a_1$, and such that $a_i=0$ for all $i>N$.
\end{abstract}
\maketitle

\renewcommand{\thefootnote}{\arabic{footnote})}
\renewcommand{\arraystretch}{1.2}

\section{Introduction}
In this paper we construct a bijection between the set of
configurations ${\bf a}=(a_i)_{i\in{\bf Z}_{\geq0}}$
satisfying the conditions
\begin{eqnarray}
&&a_i=0\hbox{ if }i>\hskip-5pt>0,\label{INT0}\\
&&a_i+a_{i+1}+a_{i+2}\leq k,\label{INT1}
\end{eqnarray}
and the set of rigged partitions $(\lambda,\rho )$, where
$\lambda=(\lambda_1,\ldots,\lambda_n)$ is a partition satisfying
$k\geq\lambda_1\geq\cdots\geq\lambda_n\geq1$,
and $\rho =(\rho _1,\ldots,\rho _n)\in{\bf Z}_{\geq0}$ is
a set of integers satisfying
\begin{equation}\label{INT2}
\rho _i\geq \rho _{i+1}\hbox{ if }\lambda_i=\lambda_{i+1}.
\end{equation}
The set of integers $\rho $ is called a rigging of the partition $\lambda$.

The bijection preserves degrees, where the degree of a configuration
${\bf a}$ is given by
\begin{equation}
d({\bf a})=\sum_{i=0}^\infty ia_i\label{DEGREE}
\end{equation}
and the degree of a rigged partition
$(\lambda,\rho)$ is given by
\begin{equation}\label{DEG}
\sum_{1\leq i<j\leq n}A_{\lambda_i,\lambda_j}+\sum_{i=1}^n\rho_i\hbox{ where }
A_{l,l'}=2{\rm min}(l,l')+{\rm max}(l+l'-k,0).
\end{equation}
Using ${\bf m}=(m_1,\ldots,m_k)$, $m_l=\sharp\{i;\lambda_i=l\}$,
one can write
\begin{eqnarray*}
Q({\bf m})&=&\sum_{1\leq i<j\leq n}A_{\lambda_i,\lambda_j}\\
&=&\frac12(A{\bf m},{\bf m})
-\frac12\sum_{l=1}^kA_{l.l}m_l.
\end{eqnarray*}
The sum over the riggings is easy because we have
\[
\sum_{\rho_1\geq\cdots\geq\rho_m\geq0}q^{\rho_1+\cdots+\rho_m}=\frac1{(q)_m}.
\]
Therefore, the bijection implies the combinatorial identity,
\begin{equation}\label{GORDON}
\sum_{\bf a}q^{d({\bf a})}=\sum_{m_1,\ldots,m_k=0}^\infty
\frac{q^{Q({\bf m})}}{\prod_{l=1}^k(q)_{m_l}},
\end{equation}
where the summation over ${\bf a}$ is under the conditions (\ref{INT0})
and (\ref{INT1}).

We also determine the image of the following two kinds of subsets
by the bijection:

\medskip
the configurations satisfying
\begin{equation}\label{SUB1}
a_0=a,\quad a_1=b;
\end{equation}

the configurations satisfying
\begin{equation}\label{SUB2}
a_i=0\hbox{ for all }i>N.
\end{equation}

We denote by $R(r_1,\ldots,r_k)$ the set of rigged partitions satisfying
\begin{equation}\label{SUB3}
\rho_i\geq r_{\lambda_i}\hbox{ for all }1\leq i\leq n.
\end{equation}
In particular, for $a,b\geq-1$ and $a+b\leq k$ we set
\begin{equation}\label{[AB]}
R[a,b]=
\begin{cases}
R(\underbrace{0,\ldots,0}_a,\underbrace{1,\ldots,b}_b,
\underbrace{b+2,\ldots,2k-2a-b}_{k-a-b})
&\hbox{if }a,b>0;\\
\emptyset&\hbox{if }a=-1\hbox{ or }b=-1.
\end{cases}
\end{equation}
The subset corresponding to (\ref{SUB1}) is given by
\begin{equation}\label{ANS}
R[a,b]\backslash(R[a-1,b+2]\cup R[a,b-1]),
\end{equation}
where $R[a-1,b+2]=R[a-1,k-a+1]$ for $a+b=k$ is understood.
The rigged partitions corresponding to (\ref{SUB2}) are characterized by
\begin{equation}\label{SUB4}
\rho_i\leq \lambda_iN-\sum_{j\not=i}A_{\lambda_i,\lambda_j}.
\end{equation}

The character of the set of rigged partitions
restricted by (\ref{SUB3}) and (\ref{SUB4}) is given by
\begin{equation}\label{FERM}
\sum_{m_1,\ldots,m_k=0}^\infty
q^{Q({\bf m})+\sum_{i=1}^kr_im_i}\prod_{1\leq l\leq k\atop m_l\not=0}
\left[lN-\sum_{i=1}^kA_{l,i}m_i+A_{l,l}-r_l+m_l\atop m_l\right].
\end{equation}
Here $\left[m\atop n\right]$ is the $q$ binomial coefficient
\[
\left[m\atop n\right]=
\begin{cases}
\prod_{i=1}^n\frac{(1-q^{m-n+i})}{1-q^i}&\hbox{if }0\leq n\leq m;\\
0&\hbox{otherwise}.
\end{cases}
\]

We denote the character corresponding to the subset $R[a,b]$
and the restriction (\ref{SUB2})
by $\chi^{(k)}_{a,b}[N]$.

In conclusion, the bijections give the following polynomial identities.
\begin{equation}\label{POL}
\sum_{\bf a}q^{d({\bf a})}
=
\begin{cases}
\chi^{(k)}_{a,b}[N]-\chi^{(k)}_{a-1,b+2}[N]-\chi^{(k)}_{a,b-1}[N]
+\chi^{(k)}_{a-1,b+1}[N]&\hbox{if }a,b>0;\\
\chi^{(k)}_{a,0}[N]-\chi^{(k)}_{a-1,2}[N]&\hbox{if }a>0,b=0;\\
\chi^{(k)}_{0,b}[N]-\chi^{(k)}_{0,b-1}[N]&\hbox{if }a=0,b>0;\\
\chi^{(k)}_{0,0}[N]&\hbox{if }a=b=0.
\end{cases}
\end{equation}
where the summation over ${\bf a}$ is under the conditions
(\ref{INT0}), (\ref{INT1}), (\ref{SUB1}) and (\ref{SUB2}).

In general, for $r\geq1$, a configuration
${\bf a}$ is called a $(k,r)$-configuration if it satisfies
\begin{equation}\label{RC}
a_i+\cdots+a_{i+r-1}\leq k.
\end{equation}
Let us discuss some physical background for this. We can think of
$a_i$ as the number of particles in the energy level $i$. If $k=r=1$,
the restriction (\ref{RC}) can be considered as Pauli's exclusion principle.
The case $(k,r)=(1,2)$ appeared in \cite{B} in the study of
the hard hexagon model in statistical mechanics on the two-dimensional lattice.
By the corner transfer matrix method,
the computation of the one point functions for
the two-dimensional lattice model reduced to the computation of
the one-dimensional configuration sums with the condition $(k,r)=(1,2)$
in (\ref{INT1}). The case of general value of $k$ (with $r=2$) appeared
in \cite{ABF}.

In representation theory, the $(k,r)$-configurations appeared
in \cite{P} as labels parametrizing a set of monomial basis
in the level $k$ irreducible highest weight representations
of the affine Lie algebras $\widehat{\mathfrak{sl}}_r$.
Very recently, a connection to Macdonald's polynomials was found \cite{FJMM},

In \cite{ABF} and also in \cite{P},
$(k,r)$-configurations are used as labels of basis of certain
infinite dimensional graded vector spaces. The grading is given by
(\ref{DEGREE}). The statistical sum (\ref{GORDON})
$\sum_{\bf a}q^{d({\bf a})}$ gives the character of these spaces.
If $r=2$, by changing slightly the definition of configurations
(this is not essential), we have the identity
\begin{equation}\label{GOR}
\sum_{\{a_i\}_{i\geq1},a_i\geq0\atop a_i+a_{i+1}\leq k}q^{\sum_{i\geq1}ia_i}
=\sum_{m_1,\ldots,m_k=0}^\infty
\frac{q^{\frac12(G{\bf m},{\bf m})}}{\prod_{l=1}^k(q)_{m_l}},
\end{equation}
where
\begin{equation}
\label{G}
G_{l,l'}=2\,{\rm min}(l,l').
\end{equation}
This is the sum side of Gordon's generalization of Roger-Ramanujan identities
(see Theorem 7.5 of \cite{A}, the case $i=k$).

In \cite{KKMM}, similar formulas for the characters
in conformal field theory are studied extensively.
The Gordon type formulas are called fermionic formulas, and formulas
in the other side of the corresponding identities
are called bosonic formulas.
In this paper, we give a fermionic formula for the $(k,3)$-configurations.
In \cite{FJMMT}, we give a different fermionic formula for
the $(k,3)$-configurations. The fermionic formulas for the general
$(k,r)$-configurations are not known. On the other hand, a bosonic formula
for the general $(k,r)$-configurations is given in \cite{FJLMM}. 

Our method for computing the one-dimensional configuration sum
is to construct a bijection between configurations ${\bf a}$
and rigged partitions $(\lambda,\rho)$. The notion of rigged configurations,
i.e., a sequence of partitions with riggings,
was introduced by \cite{KKR} in the study of Bethe Ansatz.
In this paper, we consider a single partition $\lambda$ with rigging
$\rho$. We use the term `rigged partition' for this reason.

Let us explain the meaning of rigged partitions for
one-dimensional configurations. As we have explained,
the physical interpretation of $a_i$ is
the number of particles in the energy level $i$. Our bijection gives
another way of describing a configuration as a union of particles.
Let us call them as quasi-particles in distinction with
particles in the first interpretation. In Section 2 and after,
we simply use the term `particle' since we discuss only the second
interpretation.

If $k=1$, particles and quasi-particles are the same.
In this case, if $r\geq2$, the condition (\ref{INT1}) can be understood
as a repulsive interaction between particles: two particles
cannot occupy two energy levels which are closer than $r$.
Namely, the interaction between the particles is a two-body interaction.
The lowest energy in the $m$-particle sector is
given by $rm(m-1)/2$, and the fermionic formula read as
\[
\sum_{m=1}^\infty\frac{q^{\frac r2m(m-1)}}{(q)_m}.
\]
For $k\geq2$, we introduce quasi-particles. The condition (\ref{INT1})
means no $k+1$ particles occupy energy levels in an interval of width $r$.
This is a $(k+1)$-body interaction. However, for $r=2,3$, by introducing $k$
kinds of quasi-particles, we can reform it to a two-body interaction
between the quasi-particles. We construct a bijection between
the $(k,r)$-admissible configurations $(r=2,3)$,
and the set of rigged partitions $(\lambda,\rho)$. In the sector where
the quasi-particle content is given by $\lambda=(\lambda_1,\ldots,\lambda_n)$
$(k\geq \lambda_1\geq\cdots\geq\lambda_n\geq1)$, the lowest energy is given by
$\sum_{1\leq i<j\leq k}G_{\lambda_i,\lambda_j}$ where $G_{l,l'}$ is given by
(\ref{G}) for $r=2$, or $G_{l,l'}$ replaced by $A_{l,l'}$ for $r=3$.
These are two-body interactions. In fact, it is even more.
If we renormalize the energy in each sector, i.e., if we
subtract the lowest energy, the sum over $\rho$ is the same
as in the case of free bosons.
In this way, we can reduce the system of single kind of particles
with the $(k+1)$-body interaction to the system of $k$ kinds of free particles.

Let us consider the case $r=2$.
The lowest energy configuration in the $2$-particle sector is
\begin{eqnarray*}
&&a_0,a_1,a_2,a_3,\ldots\\
&&2,\hskip6pt0,\hskip6pt0,\hskip6pt0,\ldots.
\end{eqnarray*}
We consider this as a weight $2$ quasi-particle of energy $0$.
We increase the energy of this quasi-particle one by one as follows.
\begin{eqnarray*}
&&1,\hskip6pt1,\hskip6pt0,\hskip6pt0,\ldots,\\
&&0,\hskip6pt2,\hskip6pt0,\hskip6pt0,\ldots,\\
&&0,\hskip6pt1,\hskip6pt1,\hskip6pt0,\ldots,\\
&&0,\hskip6pt0,\hskip6pt2,\hskip6pt0,\ldots.
\end{eqnarray*}
Similarly, we can define a configuration corresponding to single quasi-particle
of weight $l$ with energy $d$. We associate it with a rigged partition
$(\lambda,\rho)$ such that $\lambda=(l)$ and $\rho_1=d$.
For example, for $k\geq3$, the $(k,2)$-admissible configuration
\begin{equation}
0,\hskip6pt0,\hskip6pt2,\hskip6pt1,\hskip6pt0,\ldots
\end{equation}
corresponds to the rigged partition with $\lambda=(3)$ and $\rho_1=7$.
Because of the condition $a_i+a_{i+1}\leq k$, the weight of a quasi-particle
is at most $k$.

In general, we define the quasi-particle content of a $(k,2)$-admissible
configuration ${\bf a}$ as follows. Set $l={\rm max}(a_i+a_{i+1})$,
and let $i_1$ be the largest integer such that $a_{i_1}+a_{i_1+1}=l$.
We consider that the configuration contains a weight $l$ quasi-particle at
$(i_1,i_1+1)$. The configuration given by
\[
(M_+{\bf a})_i=
\begin{cases}
a_{i_1}-1\hbox{ if }i=i_1,\\
a_{i_1+1}+1\hbox{ if }i=i_1+1,\\
a_i\hbox{ otherwise}.
\end{cases}
\]
is also $(k,2)$-admissible and satisfies
${\rm max}\left((M_+{\bf a})_i+(M_+{\bf a})_{i+1})\right)=l$.
We call the mapping $M_+$ the right move. For example,
\begin{eqnarray*}
3,\hskip6pt0,\hskip6pt1,\hskip6pt1,\hskip6pt0,\ldots,\\
2,\hskip6pt1,\hskip6pt1,\hskip6pt1,\hskip6pt0,\ldots,\\
1,\hskip6pt2,\hskip6pt1,\hskip6pt1,\hskip6pt0,\ldots,\\
1,\hskip6pt1,\hskip6pt2,\hskip6pt1,\hskip6pt0,\ldots,\\
1,\hskip6pt1,\hskip6pt1,\hskip6pt2,\hskip6pt0,\ldots,\\
1,\hskip6pt1,\hskip6pt0,\hskip6pt3,\hskip6pt0,\ldots.
\end{eqnarray*}
In this example, the right move increases the energy of the weight $3$
quasi-particle. At first, this particle has the energy $0$.
After $5$ steps, its energy increased to $9$. We observe an acceleration
of the increment of the energy: $9-5=4$. This is equal to the energy shift
of the weight $2$ particle: at first the energy is $5$, and after
the heavy particle passes, it decreases to $1$.
In general, the energy shift when a weight $l$-particle passes a weight
$l'$-particle $(l'<l)$, is given by $G_{l,l'}$.

Let us define the quasi-particle content $\lambda=(\lambda_1,\ldots,\lambda_n)$
and the corresponding energies $\rho=(\rho_1,\ldots,\rho_n)$ of the
configuration ${\bf a}$ inductively as follows. The integer $l$ is as above.
Suppose that after $t$ steps of right moves,
the weight $l$ particle with the highest energy
is separated from the rest of the configuration. Namely, for some $j$,
$(M_+^t{\bf a})_j+(M_+^t{\bf a})_{j+1}=l$ and $(M_+^t{\bf a})_i=0$
for all $i>j+1$. Set $d=j(M_+^t{\bf a})_j+(j+1)(M_+^t{\bf a})_{j+1}$.
This is the energy of this weight $l$ particle.
Let $(\lambda_2,\ldots,\lambda_n)$ and $(\rho_2,\ldots,\rho_n)$ be
the quasi-particle content and the corresponding energies for
the rest. Then, we set $\lambda_1=l$ and
$\rho_1=t-d-\sum_{i=2}^nG_{l,\lambda_i}$.

We have sketched the bijection proof of the identity (\ref{GOR}).
In Sections 2 and 3, we construct a similar bijection for the
$(k,3)$-configurations.
\section{Particle content and rigging}
A sequence of non-negative integers 
${\bf a}=(a_i)_{i\in{\bf Z}}$ is called a configuration.
We write $({\bf a})_i$ to denote $a_i$ in ${\bf a}$.
A configuration is called finite if $a_i=0$ except for finitely many $i$,
and positively supported if $a_i=0$ for all $i<0$.
We define the energy $E({\bf a})$ and the length $|{\bf a}|$
of a finite configuration ${\bf a}$ by
\begin{eqnarray}
E({\bf a})=\sum_iia_i,\label{ENERGY}\\
|{\bf a}|=\sum_ia_i.\label{LENGTH}
\end{eqnarray}

For integer $k,r$ such that $k\geq0$ and $r\geq1$,
a configuration ${\bf a}$ is called
$(k,r)$-admissible if the following conditions are valid for all $i$.
\[
a_i+\cdots+a_{i+r-1}\leq k.
\]

In this paper we consider the case $r=3$ where we have
\begin{equation}
a_i+a_{i+1}+a_{i+2}\leq k.\label{C2}
\end{equation}
For an integer $l$ such that $0\leq l\leq k$,
a $(k,3)$-admissible configuration ${\bf a}$ is called of maximal weight $l$
if the following conditions are valid for all $i$.
\begin{eqnarray}
&&a_i+a_{i+1}\leq l,\label{C3}\\
&&a_{i-1}+2a_i+2a_{i+1}+a_{i+2}\leq k+l.\label{C4}
\end{eqnarray}

If $2l\leq k$, the condition (\ref{C4}) follows from (\ref{C3}).

\begin{dfn}\label{DEF}
We denote by $C^{(k)}$ the set of finite and $(k,3)$-admissible configurations.
We denote by $C^{(k,l)}$ the subset of $C^{(k)}$
consisting of the configurations of maximal weight less than or equal to $l$.
\end{dfn}

We abbreviate $C^{(k,l)}\backslash C^{(k,l-1)}$ to $\uC^{(k,l)}$.
The subset of $C^{(k)}$ consisting of the positively supported configurations
is denoted by $C^{(k)}_{\rm pos}$. We set
$C^{(k,l)}_{\rm pos}=C^{(k)}_{\rm pos}\cap C^{(k,l)}$.

For ${\bf a}\in C^{(k)}$, we set
\begin{equation}\label{SLDEF}
S[j,{\bf a}]=a_j+a_{j+1},\quad L[j,{\bf a}]=a_{j-1}+2a_j+2a_{j+1}+a_{j+2}.
\end{equation}
A configuration ${\bf a}$ belongs to $\uC^{(k,l)}$ if and only if
$\S i$ or $\L i$ is valid for some $i$.
It is possible that $\S i$ and $\L i$ occur at the same time.
\begin{dfn}\label{RIGHT}
We define a mapping $M_+:\uC^{(k,l)}\rightarrow\uC^{(k,l)}$ called the right
move. Let ${\bf a}\in\uC^{(k,l)}$ and let $i_1$ be the largest integer
such that $\S{i_1}$ or $\L{i_1}$ is valid.
We say that the configuration ${\bf a}$ contains a particle of weight $l$
at the highest position $i_1$.
We define a configuration $M_+{\bf a}$ by
\begin{equation}
(M_+{\bf a})_j=\begin{cases}
a_{i_1+1}+1\hbox{ if }j=i_1+1;\\
a_{i_1}-1\hbox{ if }j=i_1;\\
a_j\hbox{ otherwise}.
\end{cases}
\label{MOVE}
\end{equation}
\end{dfn}
\begin{prop}\label{RM}
If ${\bf a}\in\uC^{(k,l)}$ then $M_+{\bf a}$ belongs to $\uC^{(k,l)}$.
We have
\begin{equation}
E(M_+{\bf a})=E({\bf a})+1,\quad |M_+{\bf a}|=|{\bf a}|.
\label{EQS}
\end{equation}
\end{prop}
\begin{proof}
For notational simplicity, we write $i=i_1$. We also set $b_j=(M_+{\bf a})_j$
for all $j$. We show that $a_i>0$ so that $b_i\geq0$.
Suppose $a_i=0$. If $\S i$, then $a_{i+1}=l$.
This is a contradiction because $\S{i+1}$ then holds.
If $L[i,{\bf a}]=k+l$, then
we have $a_{i-1}+2a_{i+1}+a_{i+2}=k+l$. Since $a_{i-1}+a_{i+1}\leq k$,
we have $a_{i+1}+a_{i+2}\geq l$. This is a contradiction.

We show that $a_{i+1}+a_{i+2}+a_{i+3}<k$ so that
$b_{i+1}+b_{i+2}+b_{i+3}\leq k$. If $a_{i+1}+a_{i+2}+a_{i+3}=k$, then
we have $a_i+a_{i+1}+a_{i+2}<l$ because $a_i+2a_{i+1}+2a_{i+2}+a_{i+3}<k+l$.
This is a contradiction because neither $\S i$ nor $\L i$ holds.

After these observations it is easy to see that ${\bf b}$ belongs to
$\uC^{(k,l)}$.
The equations (\ref{EQS}) are obvious by the definition (\ref{MOVE}).
\end{proof}

\begin{exm}
\label{EX1}
The following table shows the right moves of a
configuration ${\bf a}\in C^{(3,3)}$ given by
\[
a_i=
\begin{cases}
3\hbox{ if }i=0;\\
1\hbox{ if }i=3;\\
0\hbox{ otherwise}.
\end{cases}
\]
\[
3001\rightarrow2101\rightarrow1201\rightarrow1111\rightarrow1021\rightarrow
1012\rightarrow1003.
\]
\end{exm}

One of our goals is to define a particle content of a $(k,3)$ configuration.
In Example \ref{EX1} we can think of the particle content of ${\bf a}$
to be one particle of weight $3$ and another particle of weight $1$.
In the sequence of right moves the heavy particle passes the light particle
from the left to the right.
The position of the light particle shifts by $3$ in energy.
At the same time, the right move of
the heavy particle is accelerated by the existence of the light particle
by $3$. At the first position, the energy of the heavy particle is $0$.
After the $6$ steps, it already reaches to the energy $9$.
Since the total energy difference is equal to the number of steps,
the energy shift of the light particle and the difference
between the energy shift of the heavy particle and the number od steps,
are equal, i.e., $3$ in the above example.

\begin{prop}\label{SL}
Let ${\bf a}\in C^{(k,l)}$.
Suppose that $i_1$ is the highest position of weight $l$ particle
in ${\bf a}$. If we have $\S{i_1}$,
after several right moves the highest position
will changes to $i_1+1$ $($and we have $\S{i_1+1}$ or $\L{i_1+1})$.
If $\L{i_1}$, the highest position changes to either $i_1+1$ $($and we have
$\S{i_1+1}$ or $\L{i_1+1})$, or to $i_1+2$ $($and we have $\L{i_1+2})$.
\end{prop}
\begin{proof}
While the highest position is $i_1$, the right move is nothing but
$-1,+1$ at the $i_1$-th and the $(i_1+1)$-th column.
Therefore, when the highest position changes at ${\bf b}=M_+^t{\bf a}$,
the change is such that $S[i_1+1,{\bf b}]=l$, $L[i_1+1,{\bf b}]=k+l$
or $L[i_1+2,{\bf b}]=k+l$. The change from $S[i_1,M_+^{t-1}{\bf a}]=l$
to $L[i_1+2,{\bf b}]=k+l$ is prohibited by the following lemma.
\end{proof}
\begin{lem}\label{BASIC}
Let ${\bf a}\in C^{(k,l)}$, and suppose that
\begin{equation}\label{NOL}
a_{i+1}+2a_{i+2}+2a_{i+3}+a_{i+4}=k+l
\end{equation}
for some $i$. If we have $a_i+a_{i+1}=l$, then we have $a_{i+3}+a_{i+4}=l$.
Similarly, if $a_{i+4}+a_{i+5}=l$, then we have $a_{i+1}+a_{i+2}=l$.
\end{lem}
\begin{proof}
We prove the first statement. By symmetry, the second statement follows.

We have
$a_i+2a_{i+1}+2a_{i+2}+a_{i+3}\leq k+l$.
Since $a_i+a_{i+1}=l$, we have $a_{i+1}+2a_{i+2}+a_{i+3}\leq k$.
From (\ref{NOL}), we have $a_{i+3}+a_{i+4}\geq l$. Since
${\bf a}\in C^{(k,l)}$, we have $a_{i+3}+a_{i+4}\leq l$,
and the assertion follows.
\end{proof}

Let us formulate the particle content of a configuration in general.
We set
\begin{equation}
\label{PHASE}
A_{l,l'}=2\,{\rm min}(l,l')+(l+l'-k)_+.
\end{equation}
Here, $(x)_+={\rm max}(x,0)$. The energy shift
when a heavy particle of weight $l$ passes a light particle of weight $l'$
is equal to $A_{l,l'}$. We will clarify this statement in the below.

We say a configuration ${\bf a}\in\uC^{(k,l)}$ contains
a free particle of weight $l$ at the highest position $i$
if $\S i$ is valid with $a_i\not=0$, and $a_j=0$ for all $j\geq i+2$.
Note that the right moves of such a configuration is simple.
Namely, after several changes
$-1,+1$ at the $i$-th and the $(i+1)$-th column,
the position of the free particle changes to $i+1$, and we have
$(a_i,a_{i+1},a_{i+2})=(0,l,0)$.
Then, it changes to $i+2$ and so on, each time $l$ right moves are added.
We define the energy of the free particle to be $d=ia_i+(i+1)a_{i+1}$.

\begin{prop}\label{FREE}
Let ${\bf a}\in\uC^{(k,l)}$.
If $t$ is large enough, then
$M_+^t{\bf a}$ contains a free particle of weight $l$
at the highest position $i$ for some $i$. Suppose that $(M_+^t{\bf a})_i=c$
where $1\leq c\leq l$. By the definition the energy of this particle is
$d=ic+(i+1)(l-c)$. The difference $s=d-t$ is independent of the choice of $t$.
\end{prop}
\begin{proof}
Let $i_0$ be such that
\begin{equation}
a_{i_0-1}\not=0,\quad a_j=0\hbox{ for all }j\geq i_0.
\label{ZERO}
\end{equation}
Consider the right moves of ${\bf a}$. Since the energy of ${\bf a}$
increases by $1$ in each move, in finite steps, the condition (\ref{ZERO}) will
break down. Suppose that the breakdown happens in the move from
$M_+^{t_0}{\bf a}$ to $M_+^{t_0+1}{\bf a}$.
It happens necessarily in such a way that
$(a_{i_0-1},a_{i_0})$ changing from $(l,0)$ to $(l-1,1)$.
At this stage, the configurations contain a free particle of weight $l$
at the highest position $i_0-1$. The value $s$ is independent of $t$
because both $t$ and $d$ increases by $1$ in each step.
\end{proof}
\begin{prop}\label{FREESHIFT}
Let $1\leq l'<l\leq k$. Suppose that a configuration ${\bf a}$ is such that
for some $j_1,j_2$ where $j_1<\hskip-5pt<j_2$ we have
$a_i=0$ if $i\not=j_1,j_1+1,j_2,j_2+1$ and
\[
a_{j_1}+a_{j_1+1}=l',\quad a_{j_2}+a_{j_2+1}=l.
\]
If $t$ is sufficiently large, then $M_-^t{\bf a}$ is such that
$(M_-^t{\bf a})_i=0$ if $i\not=j_3,j_3+1,j_4,j_4+1$ where
$j_3<\hskip-5pt<j_4$ is given by
\begin{eqnarray*}
&&a_{j_3}+a_{j_3+1}=l,\quad a_{j_4}+a_{j_4+1}=l',\\
&&j_3a_{j_3}+(j_3+1)a_{j_3+1}=j_2a_{j_2}+(j_2+1)a_{j_2+1}-t-A_{l.l'},\\
&&j_4a_{j_4}+(j_4+1)a_{j_4+1}=j_1a_{j_1}+(j_1+1)a_{j_1+1}+A_{l.l'}.
\end{eqnarray*}
In particular, $j_4$ is independent of $t$.
\end{prop}
\begin{proof}
For notational simplicity we set $c=a_{j_1}$. Note that $0\leq c\leq l'$.
The original configuration has a weight $l$ particle at the highest position 
$j_2$. For small $t$ the change from $M_-^t{\bf a}$ to $M_-^{t+1}{\bf a}$
is such that the energy of this particle decreases by $1$.

Case $l+l'\leq k$.
The weight $l$ particle moves until the configuration becomes of the form
\[
\ldots,0,c,l'-c,l-l'+c,l'-c,0,\ldots.
\]
We have $0+2c+2(l'-c)+(l-l'+c)=l+l'+c<k+l$ and $(l'-c)+(l-l'+c)=l$.
The configuration further changes to
\[
\ldots,0,c,l-c,c,l'-c,0,\ldots.
\]

Case $l+l'\geq k+c$.
The weight $l$ particle moves until the configuration becomes of the form
\[
\ldots,0,c,l'-c,0,k-l'+c,l+l'-k-c,0,\ldots.
\]
We have $(l'-c)+2\cdot0+2(k-l'+c)+(l+l'-k-c)=k+l$.
The configuration further changes to
\[
\ldots,0,c,l'-c,l-l',k-l+c,l+l'-k-c,0,\ldots,
\]
where $c+2(l'-c)+2(l-l')+k-l+c=k+l$, and further to
\[
\ldots,0,c,l-c,0,k-l+c,l+l'-k-c,0,\ldots.
\]

Case $k+c>l+l'>k$.
The weight $l$ particle moves until the configuration becomes of the form
\[
\ldots,0,c,l'-c,k+c-2l',l+2l'-k-c,0,\ldots.
\]
We have $c+2(l'-c)+2(k+c-2l')+l+2l'-k-c=k+l$.
The configuration changes to
\[
\ldots,0,c,l-c,k+c-l-l',l+2l'-k-c,0,\ldots.
\]

In all cases, the last configuration has a free particle at the lowest
position, and the after is simple. We can easily check that the energy
shift is equal to $A_{l,l'}$.
\end{proof}

\begin{dfn}
Let ${\bf a}\in C^{(k)}$. We define a partition
$\lambda=(\lambda_1,\ldots,\lambda_n)$ and a set of integer
$\rho =(\rho _1,\ldots,\rho _n)$ inductively with respect to
the length of configuration. We call $\lambda$ the particle content of
${\bf a}$, and $\rho $ the rigging of $\lambda$.

The inductive procedure is as follows. Let $l$ be such that
${\bf a}\in\uC^{(k,l)}$. We set $\lambda_1=l$. Let $i_0,t_0$
and $s_1=s$ be as in Proposition \ref{FREE}.
We define a new configuration $\overline{\bf a}=(\overline a_i)$
of smaller length by
\begin{equation}
\overline a_i=
\begin{cases}
(M_+^{t_0}{\bf a})_i\hbox{ if }i\leq i_0-2;\\
0\hbox{ otherwise.}
\end{cases}
\label{REMOVE}
\end{equation}
Let $\overline\lambda=(\lambda_2,\ldots,\lambda_n)$ and 
$\overline \rho =(\rho _2,\ldots,\rho _n)$ be the particle content and
its rigging. We set
\begin{equation}
\label{RIGGING}
\rho _1=s_1-\sum_{a=2}^mA_{\lambda_1,\lambda_a}.
\end{equation}
After this procedure, we define $\lambda=(\lambda_1,\overline\lambda)$
and $\rho =(\rho _1,\overline \rho )$.
\end{dfn}

We write the particle content $\lambda$, alternatively by
${\bf m}=(m_1,\ldots,m_k)$ where $\lambda=(k^{m_k},\ldots,1^{m_1})$.
Namely, $m_l$ is the number of particles in ${\bf a}$ that are of weight $l$.

The following statement is obvious by the definition.
\begin{lem}\label{OBV}
Let ${\bf a}\in C^{(k)}$, and let $(\lambda,\rho )$ be its particle content
and the rigging. Let $(\mu,s)$ be the particle content and its rigging of
$M_+{\bf a}$. Then, we have $\mu=\lambda$ and $s_i=\rho _i+\delta_{i,1}$.
\end{lem}

Let us explain the reason for the subtraction in the definition of $\rho _1$.
Suppose that a configuration ${\bf a}\in C^{(k,l)}$ is such that
$a_i=0$ for $i<0$, and $a_0=l$. Suppose further that
$i_1$ in the definition of $M_+$ is equal to $0$. Namely,
it contains a particle of weight $l$ at the highest position $0$.
We will show that the difference between the energy shift
and the number of steps when this particle
moves to the right and becomes free, is given by
$\sum_{a=2}^mA_{\lambda_1,\lambda_a}$. Suppose that after $t$ steps
of right moves the weight $l$ particle becomes free and reaches
the energy $d$. Then, the difference
$s_1=d-t$ is equal to the above sum. Namely, we have $\rho _1=0$
by the definition. In general, we will prove that
the positivity of the rigging in this normalization is equivalent to
the positivity of its support, i.e., ${\bf a}\in C^{(k)}_{\rm pos}$.

By the definition it is obvious that the particle content $\lambda$
is a partition, i.e., $\lambda_i\geq\lambda_{i+1}$ for all $i$.
It is less obvious but true that the rigging satisfies the condition
\[
\rho _i\geq \rho _{i+1}\hbox{ if }\lambda_i=\lambda_{i+1}.
\]
To prove this statement (see Proposition \ref{RIG}),
we prepare a few propositions.

In Definition \ref{RIGHT}, we defined the integer $i_1$
for a configuration ${\bf a}\in\uC^{(k,l)}$, which is the position
of the highest (or first) particle of weight $l$ in ${\bf a}$.
The right move of ${\bf a}$ is nothing but to move this particle to the right.
After finite steps, this particle becomes free. Then, we have removed
this particle from the configuration to obtain
$\overline{\bf a}$ in (\ref{REMOVE}). This is equivalent to
applying the right move $M_+$ to ${\bf a}$ infinitely many times:
\[
\overline{\bf a}=M_+^\infty{\bf a}.
\]
Applying the same procedure to $\overline{\bf a}$ and so on, we obtain
$i_2,\ldots,i_{m_l}$, which are by definition the positions of
the second particle of weight $l$, and so on. This is an inductive procedure
using many steps of right moves.
However, we can define these numbers without using right moves.

Suppose we find the integer $i_1$ from the configuration ${\bf a}$
as explained above. Instead of moving the configuration, we consider
the cut-off ${\bf b}$ of ${\bf a}$ at $i_1$:
\begin{equation}
\label{CO}
b_i=
\begin{cases}
a_i\hbox{ if }i\leq i_1-1;\\
0\hbox{ otherwise.}
\end{cases}
\end{equation}
If ${\bf b}\in\uC^{(k,l)}$, we define $\overline i_2$
to be the position of the first particle of weight $l$ in ${\bf b}$.
Continuing further while the cut-off particle still belongs
to $\uC^{(k,l)}$, we can define the numbers
$\overline i_3,\ldots,\overline i_{\overline m_l}$.
Now, we state the proposition.

\begin{prop}\label{LPC}
We follow the above setting. Then, we have the equalities
\[
\overline m_l=m_l,\quad\overline i_a=i_a(2\leq a\leq m_l).
\]
\end{prop}
\begin{proof}
It is enough to show that the position of the second particle
of weight $l$ is invariant by the right move. Let us prove this statement.
Let ${\bf a}\in C^{(k,l)}$ is such that the first particle of weight $l$
is at the position $i$ (i.e., $i_1=i$). If the position of the first
particle of weight $l$ does not change for $M_+{\bf a}$, the statement is
clear. We have two other cases, either the position of the first particle
changes to $i+1$ or to $i+2$.
In the former, the cut-off configuration obtained from $M_+{\bf a}$
is of the form
\[
\cdots,a_{i-2},a_{i-1},a_i-1,0,0,\cdots.
\]
Since ${\bf a}\in C^{(k,l)}$, this configuration satisfy
neither $S_j$ ($j\geq i-1$) nor $L_j$ ($j\geq i-2$).
Therefore, the highest position of the weight $l$ particle
is the same as that of the cut-off configuration obtained from ${\bf a}$.
In the latter, we have $L_{i+2}$ for $M_+{\bf a}$, and by Proposition \ref{SL},
we have $a_i+a_{i+1}<l$.
The cut-off configuration obtained from $M_+{\bf a}$ is of the form
\[
\cdots,a_{i-2},a_{i-1},a_i-1,a_{i+1}+1,0,\cdots.
\]
It is again easy to check that neither $S_j$ ($j\geq i-1$)
nor $L_j$ ($j\geq i-2$) is valid for this configuration.
\end{proof}
For a configuration ${\bf a}\in C^{(k,l)}$, we can thus define
the number of the weight $l$ particles $m_l$, and their positions
$i_1,\ldots,i_{m_l}$. We denote by $C^{(k,l)}_m$
the set of configurations in $C^{(k,l)}$ such that $m_l=m$.

We use the following lemma in the proof of Proposition \ref{3.1}.
\begin{lem}\label{LEM3.1}
Suppose that ${\bf a}\in C^{(k,l)}_m$.
Let $i_1,\ldots,i_m$ be the positions of the weight $l$ particles in ${\bf a}$.
Let ${\bf c}$ be the configuration obtained from ${\bf a}$
by the cut-off from the left at the column $i_m+1$. Namely,
\[
({\bf c})_i=
\begin{cases}
a_i\hbox{ if }i\geq i_m+2;\\
0\hbox{ otherwise}.
\end{cases}
\]
Then, the number of the weight $l$ particle in ${\bf c}$ is $m-1$,
and their positions are $i_1,\ldots,i_{m-1}$.
\end{lem}
\begin{proof}
In Proposition \ref{LPC} we have shown that the number and the positions of the
weight $l$ particles are determined by the cut-off procedure.
The claim of this lemma is that the cut-off from the left
in the definition of ${\bf c}$ does not affect
this procedure until we locate the $(m-1)$-th weight $l$ particle in ${\bf c}$.
To prove this it is enough to show that if
$L[i_{m-1},{\bf a}]=k+l$ then $L[i_{m-1},{\bf c}]=k+l$. This is clear if
$i_m+1<i_{m-1}-1$. Otherwise, we have $i_m+2=i_{m-1}$, and therefore,
$S[i_m,{\bf a}]=l$ and $L[i_m+2,{\bf a}]=k+l$. By Lemma \ref{BASIC}
this imples $S[i_{m-1}+1,{\bf a}]=l$. This is a contradiction.
\end{proof}

The right move $M_+$ moves the first particle
which is located at the position $i_1$. The change is
$(a_{i_1},a_{i_1+1})\rightarrow(a_{i_1}-1,a_{i_1+1}+1)$.
The number of the weight $l$ particles is invariant by this change.
The position of the first particle is either unchanged
or moves to $i_1+1$ or $i_1+2$. The positions of the other
particles are unchanged. It is natural to think of moves
of other particles of weight $l$. We want to define
$M^{(c)}_+$ ($1\leq c\leq m_l$), which changes $(a_{i_c},a_{i_c+1})$
to $(a_{i_c}-1,a_{i_c+1}+1)$. However, this is not always possible
because this change may break down the condition ${\bf a}\in C^{(k,l)}$.
In Proposition \ref{RM}, we proved that for $c=1$
the condition ${\bf a}\in C^{(k,l)}$ is preserved.
The following proposition gives an alternative answer for the case $c\geq2$.
\begin{lem}\label{WELL}
Let ${\bf a}\in C^{(k,l)}_m$. For $2\leq c\leq m$
the configuration
\[
{\bf a}^{(c)}
=M^{(c)}_+M^{(c-1)}_+\cdots M^{(1)}_+{\bf a}
\]
belongs to $C^{(k,l)}_m$.
\end{lem}
\begin{proof}
Suppose that ${\bf a}^{(c-1)}$ belongs to $C^{(k,l)}$. Set
${\bf a}^{(c-1)}=(b_i)_{i\in{\bf Z}}$. By Proposition \ref{LPC},
the position of the $c$-th particle of weight $l$ for ${\bf b}$
is equal to $i_c$, i.e., the same as ${\bf a}$.
We want to show that the change of $(b_{i_c},b_{i_c+1})$
to $(b_{i_c}-1,b_{i_c+1}+1)$ does not break the conditions
(\ref{C3}) and (\ref{C4}). Since the argument is the same for all $c\geq2$,
let us consider the case $c=2$. For simplicity we write $i=i_1$ and $i'=i_2$.
We have ${\bf a},{\bf b}=M_+{\bf a}\in C^{(k,l)}$
where
\[
b_j=
\begin{cases}
a_{i+1}+1\hbox{ if }j=i+1,\\
a_i-1\hbox{ if }j=i,\\
a_j\hbox{otherwise}.
\end{cases}
\]
We set
\[
c_j=
\begin{cases}
b_{i'+1}+1\hbox{ if }j=i'+1,\\
b_i'-1\hbox{ if }j=i',\\
b_j\hbox{otherwise}.
\end{cases}
\]
We must show that
\[
S[j,{\bf c}]\leq l,\quad L[j,{\bf c}]\leq k+l.
\]
First consider $S[j,{\bf c}]$. Since ${\bf b}\in C^{(k,l)}$,
we have to consider only the case $j=i'+1$, where
$S[j,{\bf c}]=S[j,{\bf b}]+1$. If $i'+2<i$, the positions
$(i'+1,i'+2)$ used in $S[i'+1,{\bf c}]$ are below the cut-off point
of ${\bf a}$ (see (\ref{CO}))
in the definition of the position of the second particle.
Therefore, by Proposition \ref{RM}, we have $S[j,{\bf c}]\leq l$.
In this way, the remaining case is $i'=i-2$ and $j=i-1$.
In this case, we have $S[j,{\bf c}]=S[j,{\bf a}]$ and the assertion follows.

Next consider $L[j,{\bf c}]$. The cases $j=i'+1,i'+2$ are in question
since we have $L[j,{\bf c}]=L[j,{\bf b}]+1$ for them.
Again, if the positions $(j-1,j,j+1,j+2)$ are below the cut-off point,
i.e., if $j+2<i$ we have $L[j,{\bf c}]\leq k+l$ by using
Proposition \ref{RM}. The remaining cases are
$(i',j)=(i-2,i-1),(i-2,i),(i-3,i-2),(i-3,i-1)$.
In the first two cases, we have $a_{i-2}+a_{i-1}=l$
since otherwise we must have the condition for the weight $l$ particle
in the form $a_{i-3}+2a_{i-2}+2a_{i-1}+a_i=k+l$, but the position $i$
is above the cut-off point. One can check that
except for the second case, we have $L[j,{\bf c}]=L[j,{\bf a}]$
and therefore, the assertion follows in these cases.
Finally, suppose that $L[j,{\bf c}]=k+l+1$ in the second case.
It implies that $L[i,{\bf a}]=k+l$. Recall that $a_{i-2}+a_{i-1}=l$.
By Lemma \ref{BASIC} we have $a_{i+1}+a_{i+2}=l$. This is a contradiction
because we assumed that the first particle in ${\bf a}$ is at $i$.

The invariance of $m_l$ follows from Proposition \ref{LPC}.
\end{proof}
\begin{lem}\label{CONN}
Suppose that ${\bf a}\in C^{(k,l)}_m$. For $2\leq c\leq m$ and for all
$s\geq1$, the mapping $(M^{(c)}_+)^s\cdots(M^{(1)}_+)^s{\bf a}$
is well-defined on $C^{(k,l)}_m$, and we have the equality
\[
(M^{(c)}_+)^s\cdots(M^{(1)}_+)^s{\bf a}=
(M^{(c)}_+\cdots M^{(1)}_+)^s{\bf a}.
\]
\end{lem}
\begin{proof}
The well-definedness for $s=1$ is proved in Lemma \ref{WELL}.
Set $A=M^{(c)}_+$ and $B=M^{(c-1)}_+\cdots M^{(1)}_+$.
We will show that $(AB)^s=A^sB^s$. Then, the statement of
the lemma follows by induction. It is enough to show that $AB=BA$
on the image $BC^{(k,l)}_m$ since the assertion is obtained by repeated
use of this commutativity. Let ${\bf a}\in C^{(k,l)}_m$.
By Proposition \ref{LPC}
the position of the $c$-th weight $l$ particle is the same
for ${\bf a},B{\bf a},B^2{\bf a}$. The positions of the first $c-1$ weight $l$
particles are the same for $B{\bf a}$ and $AB{\bf a}$ because the
change caused by $A$ does not alter the configuration in the region where
the first $c-1$ weight $l$ particles exist.
Therefore, the change from $B{\bf a}$ to $AB^2{\bf a}$ and that
from $B{\bf a}$ to $BAB{\bf a}$ are the same.
\end{proof}
Recall the cut-off procedure to determine the positions of the
weight $l$ particles for a configuration ${\bf a}\in C^{(k,l)}_m$.
We say the weight $l$ particles in ${\bf a}$ are free
if at each step of the cut-off procedure
we find the highest weight $l$ particle is free.
\begin{prop}\label{RIG}
Let $(\lambda,\rho )$ be the particle content and its rigging
of a configuration ${\bf a}\in C^{(k)}$. The rigging satisfies the condition
\begin{equation}\label{RL}
\rho _i\geq \rho _{i+1}\hbox{ if }\lambda_i=\lambda_{i+1}.
\end{equation}
\end{prop}
\begin{proof}
Suppose that ${\bf a}\in C^{(k,l)}_m$. By Lemma \ref{CONN}
\[
{\bf a}[t]=(M^{(m)}_+)^t\cdots(M^{(1)}_+)^t{\bf a}
\]
belongs to $C^{(k,l)}$. By Proposition \ref{FREE}, if $t$ is large enough,
in ${\bf a}[t]$ the weight $l$ particles are free. Let $d_1,\ldots,d_m$
be their energies. By the definition
$\rho _i-\rho_{i+1}=d_i-d_{i+1}-A_{l,l}$ for all $i$.
Therefore, the assertion follows from the following lemma.
\end{proof}
\begin{lem}\label{RIGCHK}
Suppose that ${\bf a}\in C^{(k,l)}$. Suppose that
for $i,i'$ such that $i\geq i'+2$ we have
$a_i+a_{i+1}=l$, $a_{i'}+a_{i'+1}=l$ and $a_j=0$ if $j\not=i,i+1,i',i'+1$.
Set $d=ia_i+(i+1)a_{i+1}$ and $d'=i'a_{i'}+(i'+1)a_{i'+1}$.
Then we have
\[
d-d'\geq A_{l,l}.
\]
\end{lem}
\begin{proof}
Recall that $A_{l,l}=2l+(2l-k)_+$. Set $a_i=a\geq0$ and $a_{i'}=b\geq0$.
We consider three cases.

Case 1: $i'\leq i-4$. We have $d-d'\geq4l+b-a$.
Since $a\leq l\leq k$, we have $4l+b-a\geq 2l+(2l-k)_+$.

Case 2: $i'=i-3$. We have $d-d'=3l+b-a$.
We have $3l+b-a\geq 2l$ because $a\leq l$.
We use the condition (\ref{C4}) for the sequence
$l-b,0,a,l-a$. It gives $3l+b-a\geq 4l-k$. The assertion follows from these
two inequalities.

Case 3: $i'=i-2$. We have $d-d'=2l+b-a$.
We have the sequence $b,l-b,a,l-a$.
Since $l-b+a\leq l$ we have $2l+b-a\geq 2l$.
The condition (\ref{C4}) gives $2l+b-a\geq 4l-k$.
The assertion follows from these two inequalities.
\end{proof}
Conversely, we have
\begin{lem}\label{CONVERSE}
Let $d_i$ $(1\leq i\leq m)$ be a set of integers
satisfying $d_i-d_{i+1}\geq A_{l,l}$ for $1\leq i\leq m-1$.
We choose $j_i,c_i$ $(1\leq i\leq m)$ such that $1\leq c_i\leq l$ and
\[
j_ic_i+(j_i+1)(l-c_i)=d_i.
\]
Then the $2m$ integers $j_i,j_i+1$ $(1\leq i\leq m)$ are distinct.
We define a configuration ${\bf a}_{\rm free}(d_1,\ldots,d_m)$ by
\[
({\bf a}_{\rm free}(d_1,\ldots,d_m))_j=
\begin{cases}
c_i\hbox{ if }j=j_i;\\
l-c_i\hbox{ if }j=j_i+1;\\
0\hbox{ otherwise.}
\end{cases}
\]
Then, the configuration ${\bf a}_{\rm free}(d_1,\ldots,d_m)$
belongs to $C^{(k,l)}_m$.
\end{lem}
\begin{proof}
First we prove that $j_i,j_i+1$ ($1\leq i\leq m$) are distinct. It is enough
to show that $j_{i+1}+1<j_i$. Without loss of generality we assume that $i=1$.
Note that $A_{l,l}\geq 2l>0$. Therefore, $j_2<j_1$.
Let us show that $j_2+1<j_1$. Suppose that $j_2+1=j_1$,
then we have $A_{l,l}\leq d_1-d_2=c_2+l-c_1<2l$. This is a contrdiction.

Set ${\bf a}={\bf a}_{\rm free}(d_1,\ldots,d_m)$
We must check the inequalities $S[i,{\bf a}]=l$ and $L[i,{\bf a}]=k+l$.
The possible cases where these inequalities are broken are
\begin{eqnarray*}
&&S[j_1-1,{\bf a}]=(l-c_2)+c_1>l,\\
&&L[j_1-1,{\bf a}]=c_2+2(l-c_2)+2c_1+(l-c_1)>k+l,\\
&&L[j_1-2,{\bf a}]=l-c_3+2c_2+2(c_2-l)+c_1>k+l,\\
&&L[j_1-1,{\bf a}]=l-c_2+2\cdot0+2c_1+(l-c_1)>k+l,\\
&&L[j_1-2,{\bf a}]=c_2+2(l-c_2)+2\cdot0+c_1>k+l,\\
&&L[j_1-2,{\bf a}]=0+2c_2+2(c_2-l)+c_1>k+l,\\
&&L[j_1-1,{\bf a}]=l-c_2+2c_1+2(c_1-l)+0>k+l.
\end{eqnarray*}
For notational simplicity we used the indices $i=1,2,3$ for $c_i$.

In each case, it is easy to lead to a contradiction to the assumption
$d_i-d_{i+1}\geq A_{l,l}={\rm max}(2l,4l-k)$.
\end{proof}

\section{Bijection between configurations and rigged partitions}
A pair of partition $\lambda=(\lambda_1,\ldots,\lambda_n)$ and 
its rigging $\rho =(\rho _1,\ldots,\rho _n)$ is called a rigged partition.
Here $n$ is a non-negative integer, $\lambda_i$ are integers satisfying
$\lambda_1\geq\cdots\geq\lambda_n>0$ and $\rho _i$ are integers satisfying
the condition (\ref{RL}). There is a unique element with $n=0$,
which we denote by $\emptyset$. The integer $n$ is specified for
each $\lambda$. In this sense we write $n=\ell(\lambda)$.
We denote by $R^{(l)}$ the set of rigged partitions satisfying
$\lambda_1\leq l$.
We denote by $R^{(l)}_{\rm pos}$ the subset of $R^{(l)}$ satisfying
$\rho _i\geq0$. We set $R^{(l)}_{\rm pos}=R^{(l)}\cap R^{(l)}_{\rm pos}$.
Note that there is a natural embedding
\[
R^{(k)}\supset R^{(k-1)}\supset\cdots\supset R^{(1)}\supset R^{(0)}
=\{\emptyset\}.
\]

In the previous section we defined a mapping
\begin{equation}
\iota:C^{(k,l)}\rightarrow R^{(l)}.
\end{equation}

We will show that this is a bijection.

In the definition of $(\lambda,\rho )$ for a configuration
${\bf a}\in C^{(k,l)}$,
we used right moves. We can define left moves and the related objects
similarly.
For example, the left move $M_-$ moves the weight $l$ particle in
${\bf a}\in \uC^{(k,l)}$ at the lowest position to the left.
To be precise, let $j_1$ be the smallest
integer such that (\ref{C3}) or (\ref{C4}) is valid for $i=j_1$.
We define ${\bf b}=M_-{\bf a}$ by
\[
b_i=\begin{cases}
a_{j_1}+1\hbox{ if }i=j_1;\\
a_{j_1+1}-1\hbox{ if }i=j_1+1;\\
a_i\hbox{ otherwise.}
\end{cases}
\]
We define the cut-off ${\bf c}$ of ${\bf a}\in \uC^{(k,l)}$ from the left
at $j_1+1$ by
\[
c_i=
\begin{cases}
a_i\hbox{ if }i\geq j_1+2;\\
0\hbox{ otherwise.}
\end{cases}
\]
Using the cut-off from the left, we can inductively determine
the number of the weight $l$ particles $m'_l$ and their positions
$j_1,\ldots,j_{m'_l}$. We can also define the mappings $M^{(c)}_-$
by changing $(a_{j_c},a_{j_c+1})$ to $(a_{j_c}+1,a_{j_c+1}-1)$.

\begin{prop}\label{3.1}
Suppose that ${\bf a}\in C^{(k,l)}$. Let $m$ and
$i_1({\bf a}),\ldots,i_m({\bf a})$
be the number and the positions of the weight $l$ particles in ${\bf a}$
with respect to the right move, and let $m'$ and
$j_1({\bf a}),\ldots,j_{m'}({\bf a})$
be the number and the positions of the weight $l$ particles in ${\bf a}$
with respect to the left move. We define the sets of integers
${\bf i}({\bf a})=\{i_1({\bf a}),\ldots,i_m({\bf a})\}$
and
${\bf j}({\bf a})=\{j_1({\bf a}),\ldots,j_{m'}({\bf a})\}$.
Then we have $m=m'$ and the equality of the sets
\begin{equation}\label{JI}
{\bf j}(M^{(m)}_+\cdots M^{(1)}_+{\bf a})={\bf i}({\bf a}).
\end{equation}
\end{prop}
Moreover, we have
\begin{equation}\label{-+}
M^{(m)}_-\cdots M^{(1)}_-M^{(m)}_+\cdots M^{(1)}_+{\bf a}={\bf a}.
\end{equation}
\begin{proof}
Set ${\bf b}=M^{(m)}_+\cdots M^{(1)}_+{\bf a}$.
First we show that $j_1({\bf b})=i_m({\bf a})$.
Set $i=i_m({\bf a})$. Then, we have $S[i,{\bf a}]=l$ or
$L[i,{\bf a}]=k+l$. Since $S[i,{\bf b}]=S[i,{\bf a}]$
and $L[i,{\bf b}]=L[i,{\bf a}]$, we have $S[i,{\bf b}]=l$ or
$L[i,{\bf b}]=k+l$. Therefore, in order to prove
$j_1({\bf b})=i_m({\bf a})$, it is enough to show that
$S[j,{\bf b}]<l$ and $L[j,{\bf b}]<k+l$ for $j<i$.
Since $i$ is the lowest position of the weight $l$ particle
with respect to the right move, we have
$S[j,{\bf b}]<l$ if $j\leq i-2$
and $L[j,{\bf b}]<k+l$ if $j\leq i-3$.
The remaining cases are $S[i-1,{\bf b}]<l$,
$L[i-2,{\bf b}]<k+l$ and $L[i-1,{\bf b}]<k+l$.
Since $b_i=a_i-1$ and $b_{i+1}=a_{i+1}+1$ these inequalities follow
from $S[i-1,{\bf a}]\leq l$, $L[i-2,{\bf a}]\leq k+l$
and $L[i-1,{\bf a}]\leq k+l$.

Now, we prove (\ref{JI}) by induction
on the length $|{\bf a}|$ of ${\bf a}$ given by (\ref{LENGTH}).
Then, the equality (\ref{-+}) follows by the definition of these mappings.

Let us consider the configuration ${\bf c}$:
\[
{\bf c}_i=
\begin{cases}
a_j&\hbox{if }j>i_m+2;\\
0&\hbox{otherwise}.
\end{cases}
\]
We have $|{\bf c}|<|{\bf a}|$. By Lemma \ref{LEM3.1}, the number of
the weight $l$ particles in ${\bf c}$ with respect to the right move is $m-1$,
and their positions are the same as $i_1({\bf a}),\ldots,i_{m-1}({\bf a})$.
Therefore, if we define ${\bf d}$ by
\[
{\bf d}_i=
\begin{cases}
b_j&\hbox{if }j>i_m+2;\\
0&\hbox{otherwise},
\end{cases}
\]
we have ${\bf d}=M_+^{(m-1)}\cdots M_+^{(1)}{\bf c}$.
By the definition the positions of the weight $l$ particles
in ${\bf d}$ with respect to the left move is
$j_2({\bf b}),\ldots,j_{m'}({\bf b})$.
Applying the induction hypothesis to ${\bf c}$, we obtain $m=m'$
and $\{i_1({\bf a}),\ldots,i_{m-1}({\bf a})\}=
\{j_2({\bf b}),\ldots,j_{m}({\bf b})\}$.
Noting that $j_1({\bf b})=i_m({\bf a})$, we obtain (\ref{JI}).
\end{proof}
By symmetry, we have
\begin{cor}\label{COR}
The mappings 
$M^{(m)}_+\cdots M^{(1)}_+$ and $M^{(m)}_-\cdots M^{(1)}_-$
on $C^{(k,l)}_m$ are inverse to each other.
\end{cor}

The inverse mapping to $\iota$,
\begin{equation}
\kappa:R^{(l)}\rightarrow C^{(k,l)},
\end{equation}
is defined by using the left move.

We construct $\kappa$ inductively on $l$ starting from
$\kappa(\emptyset)={\bf 0}$. Here ${\bf 0}$ is the configuration
such that $a_i=0$ for all $i$.

Suppose that $l>0$
Denote by $R^{(l)}_m$ the subset of $R^{(l)}$ satisfying the condition
that $\ell(\lambda)\geq m$ and $\lambda_1=\cdots=\lambda_m=l>\lambda_{m+1}$.
If $\ell(\lambda)=m$ we formally set $\lambda_{m+1}=0$ in this condition.
For $(\lambda,\rho )\in R^{(l)}_m$, we define $(\bar\lambda,\bar \rho )$ by
$\bar\lambda=(\lambda_{m+1},\ldots,\lambda_n)$ and
$\bar \rho =(\rho _{m+1},\ldots,\rho _n)$. We have $(\bar\lambda,\bar \rho )\in R^{(l-1)}$.
Suppose we have constructed $\kappa$ on $R^{(l-1)}$.
Set $\bar{\bf a}=\kappa(\bar\lambda,\bar \rho )\in C^{(k,l-1)}$.

We construct a configuration from $\bar{\bf a}$
by adding $m$ free particles of weight $l$ at appropriate energies.
Then, we use $(M^{(m)}_-\cdots M^{(1)}_-)^t$ to bring them to the correct
positions. In Example \ref{EX1}, the configuration $(\cdots3001\cdots)$
is mapped to the rigged partition $\lambda=(3,1)$ and $\rho =(0,0)$.
Let us consider the mapping $\kappa$ on this $(\lambda,\rho )$.
We have $\bar\lambda=1$ and $\bar \rho =0$. Therefore, we have
$\bar{\bf a}=(\cdots1000\cdots)$. We add a weight $3$ particle
at the energy $9$. We obtain $(\cdots1003\cdots)$. By applying
$(M^{(1)}_-)^6$ to this configuration, we obtain $(\cdots3001\cdots)$.

We now formulate this construction formally. Set
\begin{equation}\label{SR}
s_i=\rho _i+\sum_{j>i}A_{l,\lambda_j}\hbox{ for }1\leq i\leq m.
\end{equation}
For a sufficiently large $t$ we set $d_i=s_i+t$.
The condition (\ref{RL}) implies $d_i-d_{i+1}\geq A_{l,l}$
for $1\leq i\leq m-1$.
By Lemma \ref{CONVERSE} we can construct the configuration
${\bf a}_{\rm free}(d_1,\ldots,d_m)\in C^{(k,l)}_m$. If $t$ is large enough,
the sum ${\bf b}=\bar{\bf a}+{\bf a}_{\rm free}(d_1,\ldots,d_m)$
also belongs to $C^{(k,l)}_m$. We define
\begin{equation}\label{DEFK}
\kappa(\lambda,\rho )=(M^{(m)}_-\cdots M^{(1)}_-)^t{\bf b}.
\end{equation}
We have
\begin{prop}
The mappings $\iota$ and $\kappa$ are inverse to each other.
They give bijections between $C^{(k,l)}$ and $R^{(l)}$.
\end{prop}
\begin{proof}
We have already shown the well-definedness of these mappings.
Corollary \ref{COR} implies that they are inverse to each other.
\end{proof}
\begin{prop}
The energy and the length of a configuration ${\bf a}$
is given by the following formulas in terms of the corresponding
rigged partition $(\lambda,\rho )=\iota({\bf a})$.
\begin{eqnarray}
E({\bf a})&=&E_0(\lambda)+E_1(\rho ),\\[5pt]
&\hbox{where}&
E_0(\lambda)=\sum_{1\leq i<j\leq \ell(\lambda)}A_{\lambda_i,\lambda_j},\quad
E_1(\rho )=\sum_{1\leq i\leq \ell(\lambda)}\rho _i,\\
|{\bf a}|&=&\sum_{1\leq i\leq \ell(\lambda)}\lambda_i.
\end{eqnarray}
\end{prop}

The proof is straightforward.

Let $m_\alpha$ be the number of parts $\alpha$
in $\lambda$, i.e., $\lambda=(k^{m_k},(k-1)^{m_{k-1}},\ldots,1^{m_1})$. 
Using the sequence $m_i$ $(1\leq i\leq k)$, we can write $E_0(\lambda)$ as
\[
E_0(\lambda)=\frac12(A{\bf m},{\bf m})
-\sum_{1\leq\alpha\leq k}\frac12A_{\alpha,\alpha}m_\alpha.
\]
The identity (\ref{GORDON}) follows from this once we establish
the bijection between $C^{(k,l)}_{\rm pos}$ and $R^{(l)}_{\rm pos}$.
For the proof of the bijection,
the key fact is the following fact on the energy shift
when a heavy particle passes a configuration containing only lighter particles.

Fix $1\leq l'<l\leq k$. Let ${\bf a}$ be
a configuration in $C^{(k,l')}$. For a sufficiently large $j_1$ we
define ${\bf a}[j_1]\in C^{(k,l)}$ by
\begin{equation}\label{PUTL}
({\bf a}[j_1])_i=
\begin{cases}
l\hbox{ if }i=j_1;\\
a_i\hbox{ otherwise.}
\end{cases}
\end{equation}
If $t$ is sufficiently large, we can find a configuration
${\bf a}'=(a'_i)_{i\in{\bf Z}}\in C^{(k)}$ and integers $j_2$ and
$c$ $(1\leq c\leq l)$ such that $a'_i=0$ for $i\leq j_2+1$ and
\[
(M_-^t{\bf a}[j_1])_i=
\begin{cases}
0\hbox{ if }i<j_2;\\
c\hbox{ if }i=j_2;\\
l-c\hbox{ if }i=j_2+1;\\
a'_i\hbox{ otherwise}.
\end{cases}
\]
The configuration ${\bf a}'$ is independent of the choice of $(j_1,t)$.
We denote the mapping ${\bf a}\mapsto{\bf a}'$ by $P_l$.
We often drop $l$ when we fix it. The first statement is
\begin{prop}\label{L'}
In the above setting, we have
\[
P_l:C^{(k,l')}\rightarrow C^{(k,l')}.
\]
\end{prop}

The second statement is how the particle content and the rigging
change from ${\bf a}$ to ${\bf a}'=P{\bf a}$.
\begin{prop}\label{SHIFT} 
In the above setting, we set
$\iota({\bf a})=(\lambda,\rho )$ and $\iota({\bf a}')=(\mu,\rho')$.
Then we have
\begin{eqnarray}
&&\mu=\lambda,\\
&&\rho '_i=\rho _i+A_{l,\lambda_i}.
\end{eqnarray}
\end{prop}

Let us repeat what we assert in these propositions.
The left moves $M_-^t$ push down the weight $l$
particle from the energy $j_1l$ to $j_2c+(j_2+1)(l-c)$.
Differently speaking, the weight $l$ particle passes the configuration
${\bf a}\in C^{(k,l')}$ and change it to ${\bf a}'\in C^{(k,l')}$.
The particle content of the configuration ${\bf a}$ does not change.
The energy shift of the $j$-th particle, which has the weight $\lambda_j$,
is given by $A_{l,\lambda_j}$. The sum of these energy shifts is equal to the
difference between the number of steps $t$ and the energy shift
of the weight $l$ partcle:
\[
(j_1-j_2-1)l+c-t=\sum_{i=1}^{\ell(\lambda)}A_{l,\lambda_i}.
\]

{\it Proof of Propositions \ref{L'} and \ref{SHIFT} for $l=k$ or $l+l'\leq k$.}
\quad
The proof is easy because the change from ${\bf a}$
to ${\bf a}'$ is just a parallel shift of $3$ or $2$ columns, respectively.
Without loss of generality, we assume that $a_i=0$ if $i<1$ or $i>N$.
If $l=k$, the left moves of the configuration ${\bf a}$ proceed as follows.
\begin{eqnarray}
&&\ldots,0,a_1,\ldots,a_{N-2},a_{N-1},a_N,0,0,k,0,0,\ldots
\nonumber\\&\rightarrow&
\ldots,0,a_1,\ldots,a_{N-2},a_{N-1},a_N,0,k-a_N,a_N,0,0,\ldots
\nonumber\\&\rightarrow&
\ldots,0,a_1,\ldots,a_{N-2},a_{N-1},a_N,k-a_N-a_{N-1},a_{N-1},a_N,0,0,\ldots
\nonumber\\&\rightarrow&\cdots\nonumber\\&\rightarrow&
\ldots,0,a_1,a_2,k-a_1-a_2,a_1,a_2,\ldots
\nonumber\\&\rightarrow&
\ldots,0,a_1,k-a_1,0,a_1,a_2,\ldots
\nonumber\\&\rightarrow&
\ldots,0,k,0,0,a_1,a_2,\ldots
\nonumber
\label{COL3}
\end{eqnarray}
Note that $A_{k,j}=3j$ for all $j$, and this is consistent
with the energy shift caused by the parallel shift of $3$ columns.

If $l+l'\leq2k$, the left moves proceed as follows.
\begin{eqnarray}
&&\ldots,0,a_1,\ldots,a_{N-2},a_{N-1},a_N,0,l,0,0,\ldots
\nonumber\\&\rightarrow&
\ldots,0,a_1,\ldots,a_{N-2},a_{N-1},a_N,l-a_N,a_N,0,0,\ldots
\nonumber\\&\rightarrow&\cdots\nonumber\\&\rightarrow&
\ldots,0,a_1,l-a_1,a_1,a_2,\ldots,
\nonumber\\&\rightarrow&
\ldots,0,l,0,a_1,a_2,\ldots,
\end{eqnarray}
Note that $A_{l,j}=2j$ for all $1\leq j\leq l'$, and this is consistent
with the energy shift caused by the parallel shift of $2$ columns.
\qed

The proof for the case where $1\leq l'<l<k$ requires a lengthy calculation.
In the rest of this section, we prepare notations, and give
the main steps of the proof. The case checking is given in Appendix.

The main idea is to trace how the weight $l$ particle moves from
the right of the configuration ${\bf a}\in C^{(k,l')}$ to the left,
and changes ${\bf a}$ to ${\bf a}'=P{\bf a}$.
The totality of the configurations which interpolate between ${\bf a}$
and ${\bf a}'$ are of the form $M_-^{j}{\bf a}[j_1]$ in the notation of
(\ref{PUTL}). Here $j_1$ is sufficiently large, and $j$ can be an arbitrary
non-negative integer. In fact, the configuration depends only on
$d=lj_1-j$. Let us denote it by ${\bf a}^{(d)}$. Formally speaking,
we have ${\bf a}^{(\infty)}={\bf a}$ and ${\bf a}^{(-\infty)}={\bf a}'$.

\begin{dfn}
For each $d\in{\bf Z}$, we define the position of the weight $l$ particle
in the configuration ${\bf a}^{(d)}$ to be the integer $i=i(d,{\bf a})$
determined by the following condition:
\begin{eqnarray*}
&&\hbox{the equality $S[i,{\bf a}^{(d)}]=l$ or $L[i,{\bf a}^{(d)}]=k+l$ holds,}
\\&&\hbox{but neither $S[j,{\bf a}^{(d)}]=l$ nor $L[j,{\bf a}^{(d)}]=k+l$
holds for $j<i$.}
\end{eqnarray*}
A configuration ${\bf a}^{(d)}$ is called a node at $i$ if
\[
d={\rm min}\{d';i(d',{\bf a})=i\}.
\]
A node at $i$ is denoted by $S_i$ if $S[i,{\bf a}^{(d)}]=l$ holds,
and by $L_i$ if $L[i,{\bf a}^{(d)}]=k+l$.
The history for ${\bf a}\rightarrow{\bf a}'$
is the sequence of the nodes among the configurations ${\bf a}^{(d)}$.
Sometimes, we consider the history as a sequence of $S_i$ and $L_i$ forgetting
their contents as configurations.
\end{dfn}

The following properties are clear by the definition.

The history contains $S_i$ if $|i|$ is sufficiently large.
In general, $S_i$ and $L_i$ mix. A node can be $S_i$ and $L_i$
at the same time. It is also possible that neither $S_i$ nor $L_i$ is a node
After a node $S_i$ the history proceeds to either $S_{i-1}$ or $L_{i-1}$.
After a node $L_i$ (and when it is not $S_i$), the history proceeds
to either $S_{i-1}$, $L_{i-1}$ or $L_{i-2}$.
Suppose ${\bf a}^{(d_1)}$ is a node at $i$, and 
${\bf a}^{(d_2)}$ is the next node in the history.
Then, for all $j\not=i,i+1$, $({\bf a}^{(d)})_j$ is constant for
$d_2\leq d\leq d_1$. Moreover, for $d_2+1\leq d\leq d_1$
\[
({\bf a}^{(d-1)})_i=({\bf a}^{(d)})_i+1,\quad
({\bf a}^{(d-1)})_{i+1}=({\bf a}^{(d)})_{i+1}-1.
\]

\begin{exm} $k=4,l=3,l'=2$
Consider ${\bf a}\in C^{(4,2)}$ such that
\[
a_i=
\begin{cases}
1\hbox{ if }i=0,1,2;\\
0\hbox{ otherwise.}
\end{cases}
\]
We have $\iota({\bf a})=(\lambda,\rho )=((2,1),(1,0))$.
The history proceeds as
\begin{eqnarray*}
S_3&:&\cdots,0,1,1,1,0,3,0,\cdots\\
L_2&:&\cdots,0,1,1,1,1,2,0,\cdots\\
S_1&:&\cdots,0,1,1,2,0,2,0,\cdots\\
S_0&:&\cdots,0,1,2,1,0,2,0,\cdots\\
S_{-1}&:&\cdots,0,3,0,1,0,2,0,\cdots
\end{eqnarray*}
We obtain ${\bf a}'$ such that
$\iota({\bf a}')=(\mu,\rho ')=((2,1),(6,2))$.
Observe that the energy shifts are given by $A_{3,2}=5$ and $A_{3,1}=2$.
\end{exm}

The idea of the proof is to compare
the history for the case of ${\bf a}$ with that of $M_+{\bf a}$.
Suppose that ${\bf a}\in C^{(k,l')}_m$.
If $t$ is sufficiently large, the highest weight $l'$ particle
in $M_+^t{\bf a}$ is free, and the rest of the configuration belongs to
either $C^{(k,l')}_{m-1}$ or $C^{(k,l'-1)}$.
Therefore, we can reduce the problem to smaller $m$ or $l'$.
Repeating this reduction, we can finally reduce the problem to the case
when $l+l'\leq k$, which we have already proved.

Let us prepare another notational point. In the history, for a fixed $i$,
the value of ${\bf a}^{(d)}_i$ changes twice, in general, when the history
proceeds. In the above example, the value at the column $4$ is $0$ before the
history reaches the node $S_3$. At the node $S_3$, it changes to $3$, and at
the node $L_2$, it further changes to $2$. After that the value is unchanged.

The initial value is $a_i$ and the final value is $a'_i$. We denote by $a''_i$
the intermediate value. If the node $S_i$ (or $L_i$) follows after $S_{i+1}$
or $L_{i+1}$, it is of the form
\[
\ldots,a_{i-1},a_i,a''_{i+1},a'_{i+2},
\]

If $L_i$ follows after $L_{i+2}$, the values at the $(i+2)$-th and the
$(i+1)$-th columns change only once. In this case, the node $L_i$ is
of the form
\[
\ldots,a_{i-1},a_i,a_{i+1},a'_{i+2},\ldots.
\]

We give another example.
\begin{exm}
Let $k=5$, $l=4$ and $l'=3$.
We consider ${\bf a}$ given by
\[
a_i=\begin{cases}
1\hbox{ if }i=0,2,3;\\
2\hbox{ if }i=1;\\
0\hbox{ otherwise.}
\end{cases}
\]
We have $\iota({\bf a})=(\lambda,\rho )=((3,2),(2,1))$.
To see this we consider the right moves:
\begin{eqnarray*}
\ldots,0,1,2,1,1,0,\ldots\\
\ldots,0,1,1,2,1,0,\ldots\\
\ldots,0,1,1,0,3,0,\ldots
\end{eqnarray*}
Therefore, we have $\lambda_1=3$, $\rho _1=d-t-A_{3,2}=9-3-4=2$,
and $\lambda_2=2$, $\rho _2=1$.

The history proceeds as
\begin{eqnarray*}
S_5&:&\ldots0,1,2,1,1,0,0,4,0,\ldots\\
S_4&:&\ldots0,1,2,1,1,0,4,0,0,\ldots\\
L_3&:&\ldots0,1,2,1,1,2,2,0,0,\ldots\\
L_1&:&\ldots0,1,2,1,2,1,2,0,0,\ldots\\
S_0&:&\ldots0,1,3,0,2,1,2,0,0,\ldots\\
S_{-1}&:&\ldots0,4,0,0,2,1,2,0,0,\ldots.
\end{eqnarray*}
Note that $\iota({\bf a}')=((3,2),(10,6))$
From this, we observe that the energy shift of the weight $3$ particle
is $A_{4,3}=8$, and that of the weight $2$ particle is $A_{4,2}=5$.
\end{exm}

We start the proof of Propositions \ref{L'} and \ref{SHIFT}
for the case where $1\leq l'<l<k$ and $k<l+l'$.

The proof of Proposition \ref{L'} is a case checking on each possible history
for ${\bf a}\rightarrow P{\bf a}$.

Let us set up the cases to be checked.
Without loss of generality, we can assume that 
\begin{equation}\label{ASMP}
a_i=0\hbox{ for all }i\leq0\hbox{ and }a_1\not=0.
\end{equation}
In the below until we finish the proof of Proposition \ref{L'},
we keep this assumption.

By the definition it is obvious that
\begin{equation}\label{ONEL}
\hbox{each node in the history belongs to }C^{(k,l)}_1.
\end{equation}
Namely, the number of the weight $l$ particles is always $1$.

\begin{lem}\label{L1}
The history contains the nodes $S_i$ for all $i\leq 1$.
\end{lem}
\begin{proof}
First we prove that the history contains the node $S_1$. If not,
the abbreviated history goes through $L_1$ or $L_0\leftarrow L_2$.
The former implies $2a_1+2a''_2+a'_3=k+l$. (Here we consider the
case $L_1\leftarrow S_2$ or $L_1\leftarrow L_2$. However, the proof goes
similarly for $L_1\leftarrow L_3$.)
Since we have $a_1+a''_2+a'_3\leq k$, we have $a_1+a''_2\geq l$.
This implies $S_1$.

The latter implies $2a_1+a'_2=k+l$. Since $a_1+a'_2\leq l$,
we have $a_1=k$.  This implies $k\leq l$. This is a contradiction.

At $S_1$, the configuration is of the form
\[
\ldots,0,a_1,l-a_1,\ldots.
\]
It is now obvious that the  history contains $S_i$ for $i\leq0$.
\end{proof}

\begin{lem}\label{CSS}
The  history does not contain the sequence of nodes
$L_1\leftarrow L_3$.
Therefore, it contains one of the following.
\[
({\rm i})\quad S_2\qquad
({\rm ii})\quad L_2\leftarrow S_3\qquad
({\rm iii})\quad L_2\leftarrow L_3\qquad
({\rm iv})\quad L_2\leftarrow L_4.
\]
\end{lem}
\begin{proof}
By Lemma \ref{L1}, the history must contain the node $S_1$.
Therefore, if the  history contains the sequence
$L_1\leftarrow L_3$,
it contains the sequence $S_1\leftarrow L_3$.
By Lemma \ref{BASIC}, this implies $a'_4+a'_5=l$ in addition to $a_1+a_2=l$.
This is a contradiction to (\ref{ONEL}).
\end{proof}

We prove that if ${\bf a}\in C^{(k,l')}$ then ${\bf a}'=P{\bf a}\in C^{(k,l')}$
by induction. The induction goes on the length of ${\bf a}$.
We prepare induction steps as lemmas. Note that we give the proof
of the lemmas inside the big induction loop. Recall also that we assume
(\ref{ASMP}).
\begin{lem}\label{IND}
In the setting as above,
suppose that the history contains a node $S_j$ for some $j\geq2$
$($or $L_j$ for some $j\geq3)$.
Then, the configuration $\bar{\bf a}'=(\bar a'_i)_{i\in{\bf Z}}$ given by
\[
\bar a'_i=
\begin{cases}
a'_i\hbox{ if }i\geq j+2;\\
0\hbox{ otherwise,}
\end{cases}
\]
belongs to $C^{(k,l')}$.
%
\end{lem}
\begin{proof}
At $S_j$ we have
\[
S_j:\ldots,a_{j-1},a_j,l-a_j,a'_{j+2},\ldots.
\]
Consider a configuration $\bar{\bf a}\in C^{(k,l')}$ given by
\[
\bar a_i=
\begin{cases}
a_i\hbox{ if }i\geq j;\\
0\hbox{ otherwise.}
\end{cases}
\]
Since $a_1\not=0$, we have $|\bar{\bf a}|<|{\bf a}|$.
Therefore, by induction hypothesis, we have $P\bar{\bf a}\in C^{(k,l')}$.
The history for a weight $l$ particle passing $\bar{\bf a}$
from the right to the left, is obtained from that for ${\bf a}$
by cutting $a_i$ for $i\leq j-1$, before it proceeds
beyond $S_j$, where we have
\[
S_j:\ldots,0,a_j,l-a_j,a'_{j+2},\ldots.
\]
Therefore, $P\bar{\bf a}$ is obtained from $P{\bf a}$ by cutting $a'_i$
for $i\leq j+1$. In other words, $\bar{\bf a}'=P\bar{\bf a}$.
The statement follows from $P\bar{\bf a}\in C^{(k,l')}$.

The proof for the second statement is similar.
We have
\[
L_j:\ldots,a_{j-1},a_j,a''_{j+1},a'_{j+2},\ldots.
\]
If the  history goes as $L_j\leftarrow L_{j+2}$
we have $a_{j+1}$ in place of $a''_{j+1}$. We use the convention
$a''_{j+1}=a_{j+1}$ in that case. 
We consider a configuration $\bar{\bf a}\in C^{(k,l')}$ given by
\[
\bar a_i=
\begin{cases}
a_i\hbox{ if }i\geq j-1;\\
0\hbox{ otherwise,}
\end{cases}
\]
and apply the induction hypothesis to this configuration.
Until
\[
L_j:\ldots,0,a_{j-1},a_j,a''_{j+1},a'_{j+2},\ldots,
\]
the history is the same. Since $a_j+a''_{j+1}+a'_{j+2}\leq k$,
we have $a_{j-1}+a_j+a''_{j+1}\geq l$. Therefore, the history
proceeds to
\[
S_{j-1}:\ldots,0,a_{j-1},l-a_{j-1},a_{j-1}+a_j+a''_{j+1}-l,a'_{j+2},\ldots.
\]
As before, from this it follows that $\bar{\bf a}'$ belongs to $C^{(k,l')}$.
\end{proof}

Summarizing Lemmas \ref{CSS} and \ref{IND}, for the proof of
Proposition \ref{L'} it is enough to show the following inequalities:
\begin{eqnarray}
&&S[3,{\bf a}']\leq l'\hbox{ for (i--iv,)}\label{S3}\\[0pt]
&&L[4,{\bf a}']\leq k+l'\hbox{ for (i--iv),}\label{L4}\\[0pt]
&&S[4,{\bf a}']\leq l'\hbox{ for (ii-iv),}\label{S4}\\[0pt]
&&L[5,{\bf a}']\leq k+l'\hbox{ for (ii-iv),}\label{L5}\\[0pt]
&&S[5,{\bf a}']\leq l'\hbox{ for (iv),}\label{S5}\\[0pt]
&&L[6,{\bf a}']\leq k+l'\hbox{ for (iv).}\label{L6}
\end{eqnarray}
The case (\ref{S3}) for (i) follows from Lemma \ref{LA}.
The case (\ref{L4}) for (i) follows from Lemma \ref{LB}
The case (\ref{L5}) for (ii) follows from Lemma \ref{LB}.
The case (\ref{S3}) and (\ref{L4}) for (ii) follows from Lemma \ref{LC}.
The (\ref{S4}) for (ii) and (iii) follows from Lemma \ref{LCD}.
The cases (\ref{S3}) and (\ref{L4}) for (iii) follows from Lemma \ref{LD}.
The case (\ref{L5}) for (iii) follows from Lemma \ref{LEFG}.
The rest follow from Lemma \ref{LH}.

Proposition \ref{L'} is proved.

We show the commutativity of the mappings $P_l$ and $M_+$ on $C^{(k,l')}$.
This is a key step in the proof of Proposition \ref{SHIFT}.
\begin{prop}\label{COMM}
Suppose that ${\bf a}\in C^{(k,l')}$. Then, we have
$P_lM_+{\bf a}=M_+P_l{\bf a}$.
\end{prop}

This is obvious if $l=k$ or $l+l'\leq k$ because, in these cases,
as we have noted in the proof of Propositions \ref{L'} and \ref{SHIFT},
the mapping $P_l$ is a parallel shift. In the below, we assume that
$1\leq l'<l<k$ and $k<l+l'$.

We use induction in the proof of this proposition.
We use the length $|{\bf a}|$ as an induction parameter. If $|{\bf a}|=0$ ,
the assertion is clear.

Before going into the details, let us describe the steps in the
proof and prepare the setting.
Without loss of generality, we assume that ${\bf a}\in \uC^{(k,l')}$
and the highest position of the weight $l'$ particles in ${\bf a}$ is
$i=1$, i.e.,
\begin{equation}
S[1,{\bf a}]=l'\hbox{ or }L[1,{\bf a}]=k+l',\label{HIGH}
\end{equation}
and
\begin{equation}
\hbox{neither }S[i,{\bf a}]=l'\hbox{ nor }
L[i,{\bf a}]=k+l'\hbox{ holds for }i>1.\label{LOW}
\end{equation}

In order to know about $M_+P{\bf a}$, we need to know the highest position $i'$
of the weight $l'$ particles in ${\bf a}'=P{\bf a}$.
We show that $i'=3,4$ or $5$ depending only on the history for
${\bf a}\rightarrow{\bf a}'$. In order to know about $PM_+{\bf a}$,
we compare $(M_+{\bf a})^{(d)}$ with ${\bf a}^{(d)}$. The comparison is
not very difficult because $M_+{\bf a}$ is obtained from ${\bf a}$ by changing
$(a_1,a_2)$ to $(a_1-1,a_2+1)$. The main point is to know how the change of
$(a_1,a_2)$ to $(a_1-1,a_2+1)$ makes a difference in 
$(M_+{\bf a})^{(d)}$ compared with ${\bf a}^{(d)}$.
Two configurations are the same except at the columns $1$ and $2$,
before the node $S_i$ or $L_i$ with $i\leq 4$ appear in the history
for ${\bf a}\rightarrow{\bf a}'$.
This is because before that happens we have
\[
({\bf a}^{(d)})_j=a_j\hbox{ and }(M_+{\bf a})^{(d)}_j=(M_+{\bf a})_j
\hbox{ for }j\leq 5,
\]
and therefore, the difference at the columns $1$ and $2$ makes
no difference between $(M_+{\bf a})^{(d)}$ and ${\bf a}^{(d)}$
in the region $i\geq6$.

In the proof, we will see also that after the node $S_i$ or $L_i$
with $i\leq0$ appear in the history, $(M_+{\bf a}^{(d)})_j=({\bf a}^{(d)})_j$
for $j\leq2$. Namely, the difference at the columns $1$ and $2$
disappear. Therefore, the comparison is necessary only in the finite region
of $i$. Possible histories (considered as sequences of $S_i$ or $L_i$)
in this finite region is finite. We will check all these cases one by one.

\begin{prop}
We follow the above setting. Consider the history for
${\bf a}\rightarrow{\bf a}'$.
It does not contain $S_2$ nor $L_2\leftarrow L_3\leftarrow L_5$.
\end{prop}
We prove this proposition in Lemmas \ref{LEM0}--\ref{LEM5}.

Since the history does not contain $S_2$, it contains $L_2$
or $L_1\leftarrow L_3$. Since $L_2\leftarrow L_3\leftarrow L_5$
is also out, we have the following cases.
\begin{eqnarray*}
&&\hbox{(A) $L_2\leftarrow S_3\leftarrow S_4$ where $L_3$ is not a node,}\\
&&\hbox{(B) $L_2\leftarrow S_3\leftarrow L_4$ where $L_3$ is not a node,}\\
&&\hbox{(C) $L_2\leftarrow L_3\leftarrow S_4$,}\\
&&\hbox{(D) $L_2\leftarrow L_3\leftarrow L_4$,}\\
&&\hbox{(E) $L_2\leftarrow L_4$,}\\
&&\hbox{(F) $L_1\leftarrow L_3$,}
\end{eqnarray*}
where in all cases, $S_2$ is not a node.
\begin{rem}
In Cases (E) and (F), the history does not contain $S_2$
by the definition. In other cases, we assume that $S_2$
is not contained.
\end{rem}

{\it Case} (A).

The assumption that $L_3$ is not a node is equivalent to $a_2+a_3<k-l$
because we have $a_2+2a_3+2a''_4+a'_5<k+l$ and $a_3+a''_4=a''_4+a'_5=l$.

\begin{lem}\label{LEMD}
Suppose that the  history for ${\bf a}\rightarrow{\bf a}'$
contains the sequence
\[
L_2\leftarrow S_3\leftarrow S_4.
\]
We also assume that $L_3$ is not contained.
The  history continues to either
\[
S_1:\ldots,a_0,a_1,a''_2,a'_3,\ldots,
\]
or
\[
L_1:\ldots,a_0,a_1,a''_2,a'_3,\ldots,
\]
or
\[
L_0:\ldots,a_0,a_1,a'_2,a'_3,\ldots.
\]
For simplicity of notation let us denote $a''_2=a'_2$ in the last case.

In all cases, the  history for $M_+{\bf a}\rightarrow PM_+{\bf a}$
contains the node
\[
S_4:\ldots,a_1-1,a_2+1,a_3,a_4,a''_5,a'_6,\ldots,
\]
and it continues as
\begin{eqnarray*}
S_3&:&\ldots,a_1-1,a_2+1,a_3,a''_4,a'_5,a'_6,\ldots,\\
L_2&:&\ldots,a_1-1,a_2+1,a''_3-1,a'_4+1,a'_5,a'_6,\ldots,\\
S_1&:&\ldots,a_1-1,a''_2+1,a'_3-1,a'_4+1,a'_5,a'_6,\ldots.
\end{eqnarray*}
The mapping $M_-$ brings the last configuration to
\[
\ldots,a_1,a''_2,a'_3-1,a'_4+1,a'_5,a'_6,\ldots,
\]
and, after this, two histories are identical except for
the difference at the third and fourth columns, i.e.,
$(a'_3,a'_4)$ or $(a'_3-1,a'_4+1)$.
\end{lem}
\begin{proof}
Before $S_4$ two histories are the same.
In particular, the last node before $S_4$ is $S_5$ or $L_5$.
In both cases, the change takes place at the columns $5$ and $6$.

To see that $S_4$ appears as a node, it is enough to show that
$a_2+1+2a_3+2a_4+a''_5<k+l$. If $a_2+1+2a_3+2a_4+a''_5=k+l$,
since $a_4+a''_5=l$, we have $a_2+2a_3+a_4=k-1$.
By Lemma \ref{LEM1} we have $a_2+a_3\leq k-l$. Therefore, we have
$a_3+a_4\geq l-1\geq l'$. This is a contradiction.

Now, we assume that $a_2+a_3<k-l$. To see that $S_3$ appears as a node,
we must show that
\begin{eqnarray}
a_2+1+2a_3+2a''_4+a'_5\leq k+l,\label{EQ1}\\
a_1-1+2(a_2+1)+2a_3+a''_4<k+l.\label{EQ2}
\end{eqnarray}
Since $L_3$ is not contained in the history, we have (\ref{EQ1}).
If $a_1-1+2(a_2+1)+2a_3+a''_4=k+l$, since $a''_4+a'_5=l$,
by Lemma \ref{BASIC} we have $a_1+a_2=l$. This is a contradiction.
We have shown (\ref{EQ2}).

To see that $L_2$ appears as a node, we need
$(a_1-1)+2(a_2+1)+2(a''_3-1)+a'_4+1=k+l$,
$(a_2+1)+(a''_3-1)\leq l$ and
$a_0+2(a_1-1)+2(a_2+1)+a''_3-1<k+l$.
These are obvious.

To see that $S_1$ appears as a node, first note that by Lemma \ref{S1MUST}
we have $a_1+a''_2=l$ or $a_1+a'_2=l$. We have also
$a_0+2(a_1-1)+2(a''_2+1)+a'_3-1<k+l$ and
$a_{-1}+2a_0+2(a_1-1)+a_2+1<k+l$. Thus we have the node $S_1$.

Finally, in one step, the columns $(a_1-1,a''_2+1)$ change to $(a_1,a''_2)$,
and two histories coincide after that except for the third and the fourth 
columns.
\end{proof}
\begin{prop}\label{PROPA}
We follow the setting as given by $(\ref{HIGH})$ and $(\ref{LOW})$.
Suppose that the  history for ${\bf a}\rightarrow{\bf a}'$
contains the sequence
\[
L_2\leftarrow S_3\leftarrow S_4.
\]
Suppose also that $L_3$ does not appear as a node.
Then we have $PM_+{\bf a}=M_+P{\bf a}$.
The highest position of the weight $l'$ particles in ${\bf a}'=P{\bf a}$ is
at the column $3$.
\end{prop}
\begin{proof}
We define a configuration $\tilde{\bf a}\in C^{(k,l'-1)}$ by
a cut-off from ${\bf a}$:
\[
\tilde a_i=
\begin{cases}
a_i\hbox{ if }i\geq2;\\
0\hbox{ otherwise,}
\end{cases}
\]
and consider the history corresponding to this configuration.

We have the node
\[
L_2:\ldots,a_1,a_2,a''_3,a'_4,\ldots.
\]
Therefore, $a_1+2a_2+2a''_3+a'_4=k+l$, and
since $a_1+a_2+a''_3\leq k$ we have $a_2+a''_3+a'_4\geq l$.
The history for $\tilde{\bf a}\rightarrow P\tilde{\bf a}$ contains the node
\[
S_3:\ldots,0,a_2,a_3,a''_4,a'_5,\ldots.
\]
Since $P\tilde{\bf a}\in C^{(k,l'-1)}$,
this implies that the configuration $\tilde{\bf a}'$ given by
\[
\tilde a'_i=
\begin{cases}
a'_i\hbox{ if }i\geq5;\\
0\hbox{ otherwise},
\end{cases}
\]
belongs to $C^{(k,l'-1)}$.

Now, we show that $PM_+{\bf a}=M_+P{\bf a}$.
By Lemma \ref{LEMD} we know that $PM_+{\bf a}$
is obtained from $P{\bf a}$ by changing $(a'_3,a'_4)$
to $(a'_3-1,a'_4+1)$. On the other hand, by Lemma \ref{S1MUST}
we have $a'_3+a'_4=l'$. Therefore, to prove $PM_+{\bf a}=M_+P{\bf a}$
it is enough to show that $S[i,{\bf a}']<l'$ and $L[i,{\bf a}']<k+l'$
for $i\geq 4$. The former for $i\geq5$ and the latter for $i\geq6$
follow from $\tilde{\bf a}'\in C^{(k,l'-1)}$, and the rest follows from
$PM_+{\bf a}\in C^{(k,l')}$.
\end{proof}

{\it Case} (B).

\begin{lem}\label{LEME}
Suppose that the  history for ${\bf a}\rightarrow{\bf a}'$
contains the sequence $L_2\leftarrow S_3\leftarrow L_4$, but not $L_3$.
The  history for $M_+{\bf a}\rightarrow PM_+{\bf a}$
contains the node
\[
L_4:\ldots,a_1-1,a_2+1,a_3,a_4,a''_5,a'_6,\ldots.
\]
If $a_1+2a_2+2a_3+a''_4<k+l-1$, it continues as
\begin{eqnarray*}
S_3&:&\ldots,a_1-1,a_2+1,a_3,a''_4,a'_5,a'_6,\ldots,\\
L_2&:&\ldots,a_1-1,a_2+1,a''_3-1,a'_4+1,a'_5,a'_6,\ldots,\\
S_1&:&\ldots,a_1-1,a''_2+1,a'_3-1,a'_4+1,a'_5,a'_6,\ldots.
\end{eqnarray*}
If $a_1+2a_2+2a_3+a''_4=k+l-1$, we have $a''_3=a_3+1$ and $a'_4=a''_4-1$, and
the part of the history, $L_2\leftarrow S_3$, is replaced by only
\[
L_2:\ldots,a_1-1,a_2+1,a''_3-1,a'_4+1,a'_5,a'_6,\ldots.
\]

In both cases, the mapping $M_-$ brings the configuration $S_1$ to
\[
\ldots,a_1,a''_2,a'_3-1,a'_4+1,a'_5,a'_6,\ldots,
\]
and, after this, two histories are identical except for
the difference at the third and fourth columns, i.e.,
$(a'_3,a'_4)$ or $(a'_3-1,a'_4+1)$.
\end{lem}
\begin{proof}
To see that $L_4$ is a node, we need to show $a_2+1+2a_3+2a_4+a''_5<k+l$.
We have $a_2+2a_3+2a_4+a''_5<k+l$ because $L_4$ is a node in the history
for ${\bf a}\rightarrow{\bf a}'$. Suppose that $a_2+2a_3+2a_4+a''_5=k+l-1$.
We have $a_2+2a_3+2(a_4+1)+a''_5-1=k+l$. Since
$a_2+2a_3+2a''_4+a'_5\leq k+l$, we have $a''_4=a_4+1$. This is a contradiction
because $a_3+a_4=a_3+a''_4-1=l-1\geq l'$.

Now, we use that the history for ${\bf a}\rightarrow {\bf a}'$ has the node
$S_3$ but not $L_3$. The only obstruction for the existence of the node
$S_3$ in the history for $M_+{\bf a}\rightarrow PM_+{\bf a}$, is the value
of $a_1-1+2(a_2+1)+2a_3+a''_4=a_1+2a_2+2a_3+a''_4+1\leq k+l$.
If $a_1+2a_2+2a_3+a''_4<k+l-1$, we have $a_1-1+2(a_2+1)+2a_3+a''_4<k+l$.
Therefore, the history for $M_+{\bf a}\rightarrow PM_+{\bf a}$
has the node $S_3$. After this node the argument is
the same as in Lemma \ref{LEMD}.
If $a_1+2a_2+2a_3+a''_4=k+l-1$, the history has the node
\[
L_2:\ldots,a_1-1,a_2+1,a_3,a''_4,a'_5,a'_6,\ldots.
\]
We have $a_1+2a_2+2(a_3+1)+a''_4-1=k+l$ and $a_1+2a_2+2a''_3+a'_4+1=k+l$.
Therefore, we have $a''_3=a_3+1$ and $a'_4=a''_4-1$. The statement follows
from this observation.
\end{proof}
\begin{rem}
If $a_1+2a_2+2a_3+a''_4=k+l-1$ in Lemma \ref{LEME}, the node $L_2$
can be also written as
\[
\ldots,a_1-1,a_2+1,a_3,a''_4,a'_5,a'_6,\ldots.
\]
In other words, we can think of the history as containing the sequence
$L_2\leftarrow S_3$, where the number of steps from $S_3$ to $L_2$ is $0$.
A similar statement holds in some of other cases below. We do not repeat the
remark.
\end{rem}
\begin{prop}\label{PROPB}
We follow the setting as given by $(\ref{HIGH})$ and $(\ref{LOW})$.
Suppose that the  history for ${\bf a}\rightarrow{\bf a}'$
contains the sequence
\[
L_2\leftarrow S_3\leftarrow L_4.
\]
We also assume that $L_3$ does not appear as a node.
Then we have $PM_+{\bf a}=M_+P{\bf a}$.
The highest position of the weight $l'$ particles in ${\bf a}'=P{\bf a}$ is
at the column $3$.
\end{prop}
\begin{proof}
We use Lemmas \ref{S1MUST} and \ref{LEME}.
After the node $S_3$ the proof is the same as that of Proposition \ref{PROPA}.
\end{proof}

{\it Cases} (C) {\it and} (D).
\begin{lem}\label{LEMH}
Suppose that the  history for ${\bf a}\rightarrow{\bf a}'$
contains the sequence
\[
L_2\leftarrow L_3\leftarrow S_4(\hbox{or }L_4).
\]
We assume that $S_2$ does not appear as a node.
If the history continues as $L_0\leftarrow L_2$, we formally set $a''_2=a'_2$.
If the history contains $L_4\leftarrow L_6$, we formally set $a''_5=a_5$.

If $a_2+2a_3+2a_4+a''_5<k+l-1$,
the  history for $M_+{\bf a}\rightarrow PM_+{\bf a}$
contains the sequence
\begin{eqnarray*}
S_4\hbox{ or }L_4&:&\ldots,a_1-1,a_2+1,a_3,a_4,a''_5,a'_6,\ldots,\\
L_3&:&\ldots,a_1-1,a_2+1,a_3,a''_4-1,a'_5+1,a'_6,\ldots,\\
L_2&:&\ldots,a_1-1,a_2+1,a''_3,a'_4-1,a'_5+1,a'_6,\ldots,\\
S_1\hbox{ or }L_1&:&\ldots,a_1-1,a''_2+1,a'_3,a'_4-1,a'_5+1,a'_6,\ldots.
\end{eqnarray*}
If $a_2+2a_3+2a_4+a''_5=k+l-1$, then we have $a''_4-1=a_4$ and $a'_5+1=a''_5$,
and the part of the history, $L_3\leftarrow S_4(\hbox{or }L_4)$,
is replaced by only $L_3$.

In both cases, the mapping $M_-$ brings the configuration
$S_1(\hbox{or }L_1)$ to
\[
\ldots,a_1,a''_2,a'_3,a'_4-1,a'_5+1,a'_6,\ldots,
\]
and after this, two histories coincide except
for the fourth and the fifth columns.
\end{lem}
\begin{proof}
Since the history for ${\bf a}\rightarrow{\bf a}'$ has the node
$S_4$(or $L_4$), we have $a_2+2a_3+2a_4+a''_5\leq k+l-1$. 
If $a_2+2a_3+2a_4+a''_5\leq k+l-2$, i.e., 
$(a_2+1)+2a_3+2a_4+a''_5<k+l$, the history for
$M_+{\bf a}\rightarrow PM_+{\bf a}$ also has $S_4$(or $L_4$) as a node,
and, since $(a_1-1)+2(a_2+1)+a_3+(a''_4-1)<k+l$, it 
proceeds to $L_3$. If $a_2+2a_3+2a_4+a''_5=k+l-1$, we have
$a_2+2a_3+2(a_4+1)+(a''_5-1)=k+l$ and $(a_2+1)+2a_3+2a_4+a_5''=k+l$.
It implies that $a''_4=a_4+1$ and
$a'_5=a''_5-1$, and the history contains $L_3$ without $S_4$(or $L_4$).
In both cases, using the assumption that the  history for 
${\bf a}\rightarrow{\bf a}'$ contains $L_2$ but not $S_2$, we have
$a_0+2(a_1-1)+2(a_2+1)+a''_3<k+l$ and $(a_2+1)+a''_3\leq l$. Therefore,
the history for $M_+{\bf a}\rightarrow PM_+{\bf a}$ proceeds to $L_2$.

Next, we show that it further proceeds to $S_1$(or $L_1$),
i.e., it does not proceeds to
\[
L_0:\ldots,a_{-1},a_0,a_1-1,a''_2+2,a'_3-1,\ldots.
\]
There are three cases of the history for ${\bf a} \to {\bf a}'$:
(i) $S_1\leftarrow L_2$, (ii) $L_1\leftarrow L_2$ 
and (iii) $L_0\leftarrow L_2$. The cases (i) and (ii) is straightforward.
The case (iii) follows from Lemma \ref{LEMF}.

The rest of proof is the same as Lemma \ref{LEMD}.
\end{proof}
\begin{prop}\label{PROPC}
We follow the setting as given by $(\ref{HIGH})$ and $(\ref{LOW})$.
Suppose that the  history for ${\bf a}\rightarrow{\bf a}'$
contains the sequence
\[
L_2\leftarrow L_3\leftarrow S_4(\hbox{or }L_4).
\]
We also assume that $S_2$ does not appear as a node.
Then we have $PM_+{\bf a}=M_+P{\bf a}$.
The highest position of the weight $l'$ particles in ${\bf a}'=P{\bf a}$ is
at the column $4$.
\end{prop}
\begin{proof}
Since we have Lemmas \ref{LEMG}, \ref{LEMG'} and \ref{LEMH},
it is enough to repeat the argument in the proof of Proposition \ref{PROPA}.
\end{proof}

{\it Case} (E).

\begin{lem}\label{LEMJ}
Suppose that the  history for ${\bf a}\rightarrow{\bf a}'$
contains the sequence $L_2\leftarrow L_4$
We assume that $S_2$ does not appear as a node.
If the history continues as $L_0\leftarrow L_2$, we formally set $a''_2=a'_2$.
If the history contains $L_4\leftarrow L_6$, we formally set $a''_5=a_5$.

The  history for $M_+{\bf a}\rightarrow PM_+{\bf a}$
contains the sequence
\begin{eqnarray*}
L_4&:&\ldots,a_1-1,a_2+1,a_3,a_4,a''_5,a'_6,\ldots,\\
L_2&:&\ldots,a_1-1,a_2+1,a_3,a'_4-1,a'_5+1,a'_6,\ldots,\\
S_1\hbox{ or }L_1&:&\ldots,a_1-1,a''_2+1,a'_3,a'_4-1,a'_5+1,a'_6,\ldots.
\end{eqnarray*}
The mapping $M_-$ brings the configuration $S_1(\hbox{or }L_1)$ to
\[
\ldots,a_1,a''_2,a'_3,a'_4-1,a'_5+1,a'_6,\ldots,
\]
and after this, two histories coincide except
for the fourth and the fifth columns.
\end{lem}
\begin{proof}
We show that $(a_2+1)+2a_3+2a_4+a''_5<k+l$. If $(a_2+1)+2a_3+2a_4+a''_5=k+l$,
we have $a_2+2a_3+2(a_4+1)+(a_5''-1)=k+l$. This implies $a'_4=a_4+1$ and
$a'_5=a''_5-1$. Then, we have $a_1+2a_2+2a_3+a_4=a_1+2a_2+2a_3+a'_4-1
=k+l-1\geq k+l'$. This is a contradiction. This implies that
the history for $M_+{\bf a}\rightarrow PM_+{\bf a}$ has $L_4$ as a node.

Setting formally $a''_3=a_3$ and $a''_4=a'_4$, we can repeat the rest of
the proof of Lemma \ref{LEMH}.
\end{proof}
\begin{prop}\label{PROPE}
We follow the setting as given by $(\ref{HIGH})$ and $(\ref{LOW})$.
Suppose that the  history for ${\bf a}\rightarrow{\bf a}'$
contains the sequence
\[
L_2\leftarrow L_4.
\]
Then we have $PM_+{\bf a}=M_+P{\bf a}$.
The highest position of the weight $l'$ particles in ${\bf a}'=P{\bf a}$ is
at the column $4$.
\end{prop}
\begin{proof}
Since we have Lemmas \ref{LEMI} and \ref{LEMJ},
it is enough to repeat the argument in the proof of Proposition \ref{PROPA}.
\end{proof}

{\it Case} (F).

\begin{lem}\label{LEMN}
Suppose that the  history for ${\bf a}\rightarrow{\bf a}'$
contains the sequence
\[
L_1\leftarrow L_3\leftarrow S_4(\hbox{or }L_4).
\]
If the history contains $L_4\leftarrow L_6$, we formally set $a''_5=a_5$.

If $a_2+2a_3+2a_4+a''_5<k+l-1$,
the  history for $M_+{\bf a}\rightarrow PM_+{\bf a}$
contains the sequence
\begin{eqnarray*}
S_4\hbox{ or }L_4&:&\ldots,a_0,a_1-1,a_2+1,a_3,a_4,a''_5,a'_6,\ldots,\\
L_3&:&\ldots,a_0,a_1-1,a_2+1,a_3,a''_4-1,a'_5+1,a'_6,\ldots,\\
L_1&:&\ldots,a_0,a_1-1,a_2+1,a'_3,a'_4-1,a'_5+1,a'_6,\ldots.
\end{eqnarray*}
If $a_2+2a_3+2a_4+a''_5=k+l-1$, then we have $a''_4-1=a_4$ and $a'_5+1=a''_5$,
and the part of the history, $L_3\leftarrow S_4(\hbox{or }L_4)$,
is replaced by only $L_3$.

In both cases, the mapping $M_-$ brings the configuration
$L_1$ to
\[
\ldots,a_0,a_1,a_2,a'_3,a'_4-1,a'_5+1,a'_6,\ldots,
\]
and after this, two histories coincide except
for the fourth and the fifth columns.
\end{lem}
\begin{proof}
The proof is the same as Lemma \ref{LEMH} until
the history for $M_+{\bf a}\rightarrow PM_+{\bf a}$
reaches the node $L_3$. Now, to see that it proceeds to the node $L_1$,
it is enough to show that $a_2+a'_3<l$. The proof for this statement is the
same as Lemma \ref{LEM3} by setting $a''_3=a'_3$. After this node,
the statement of the lemma is clear.
\end{proof}
\begin{lem}\label{LEMO}
Suppose that the  history for ${\bf a}\rightarrow{\bf a}'$
contains the sequence
\[
L_1\leftarrow L_3\leftarrow L_5.
\]
If the history contains $L_5\leftarrow L_7$, we formally set $a''_6=a_6$.

The  history for $M_+{\bf a}\rightarrow PM_+{\bf a}$
contains the sequence
\begin{eqnarray*}
L_5&:&\ldots,a_0,a_1-1,a_2+1,a_3,a_4,a_5,a''_6,a'_7,\ldots,\\
L_3&:&\ldots,a_0,a_1-1,a_2+1,a_3,a_4,a'_5-1,a'_6+1,a'_7,\ldots,\\
L_1&:&\ldots,a_0,a_1-1,a_2+1,a'_3,a'_4,a'_5-1,a'_6+1,a'_7,\ldots.
\end{eqnarray*}
The mapping $M_-$ brings the configuration $L_1$ to
\[
\ldots,a_0,a_1,a_2,a'_3,a'_4,a'_5-1,a'_6+1,a'_7,\ldots,
\]
and after this, two histories coincide except
for the fifth and the sixth columns.
\end{lem}
The proof is straightforward.
\begin{prop}\label{PROPF}
We follow the setting as given by $(\ref{HIGH})$ and $(\ref{LOW})$.
Suppose that the  history for ${\bf a}\rightarrow{\bf a}'$
contains the sequence
\[
L_1\leftarrow L_3.
\]
Then we have $PM_+{\bf a}=M_+P{\bf a}$.
If the history contains the sequence
$L_3\leftarrow S_4(\hbox{or }L_4)$
$($resp., $L_3\leftarrow L_5)$,
the highest position of the weight $l'$ particles in ${\bf a}'=P{\bf a}$ is
at the column $4$ $($resp., $5)$.
\end{prop}
\begin{proof}
Since we have Lemmas \ref{LEMN} through \ref{LEMO},
it is enough to repeat the argument in the proof of Proposition \ref{PROPA}.
\end{proof}

By Propositions \ref{PROPA} through \ref{PROPF}, we have finished the
proof of Proposition \ref{COMM}.

\bigskip
{\it Proof of Proposition \ref{SHIFT}.}
We use an induction on the number of the weight $l'$ particles of ${\bf a}$.
We also use an induction on the position of the highest weight $l'$ particle:
if the assertion is valid for $M_+{\bf a}$ then it is valid for ${\bf a}$
by Lemma \ref{OBV} and Proposition \ref{COMM}. In Proposition \ref{FREE}
we have shown that the right moves on a configuration ${\bf a}$ separate
a weight $l'$ particle at the highest position. Taking this separation large
enough the mapping $P$ on $M_+^t{\bf a}$ can be separately given for the weight
$l'$ free particle, and the rest, which has less weight $l'$ particles.
Since the phase shift for a free particle
is given by Proposition \ref{FREESHIFT}, the proof is over.
\qed

\bigskip
Finally we prove
\begin{thm}
The mappings $\iota$ and $\kappa$ give the bijections between
$C^{(k,l)}_{\rm pos}$ and $R^{(l)}_{\rm pos}$.
\end{thm}
\begin{proof}
We will show
\begin{equation}\label{POSI}
\iota(C^{(k,l)}_{\rm pos})\subset R^{(l)}_{\rm pos},
\end{equation}
and
\begin{equation}\label{POSK}
\kappa(R^{(l)}_{\rm pos})\subset C^{(k,l)}_{\rm pos}.
\end{equation}

Suppose that ${\bf a}$ and $(\lambda,\rho )$ are mapped by the bijections
$\iota$ and $\kappa$ to each other.
The positivity for ${\bf a}\in C^{(k,l)}_{\rm pos}$ is that
\begin{equation}\label{POSC}
a_j=0\hbox{ for all }j<0,
\end{equation}
and the positivity for  $(\lambda,\rho )\in R^{(l)}_{\rm pos}$ is that
\begin{equation}\label{POSR}
\rho _i\geq0\hbox{ for all }1\leq i\leq \ell(\lambda).
\end{equation}

We use an induction on $l$ and the number of the weight $l$ particles,
which we denote by $m_l$.
Suppose that $m_l>1$ for ${\bf a}\in C^{(k,l)}_{\rm pos}$.
The right move $M_+$ does not change the condition (\ref{POSC}). We also
have that $\rho _i\geq \rho _{i+1}$ if $\lambda_i=\lambda_{i+1}$. Therefore,
$\iota({\bf a})\in R^{(l)_{\rm pos}}$ follows by induction.
Therefore, for the proof of (\ref{POSI}), we can assume that $m_l=1$.

The mapping $\kappa$ is defined in (\ref{DEFK}).
Note that (see Lemma \ref{CONN})
\[
(M^{(m_l)}_-)^s\cdots(M^{(1)}_-)^s{\bf a}
=(M^{(m_l)}_-\cdots M^{(1)}_-)^s{\bf a}.
\]
Suppose that $m_l>1$. By the definition if the configuration
$(M^{(1)}_-)^s{\bf a}$ is positively supported then 
$(M^{(m_l)}_-\cdots M^{(1)}_-)^s{\bf a}$ is also positively supported.
Therefore, for the proof of (\ref{POSK}), we can assume that $m_l=1$.

In the case $m_l=1$ the bijectivity follows from Proposition \ref{SHIFT}
by the following reason. It is enough to show the equivalence of the condition
\[
a_i=\begin{cases}
0(i\leq-1);\\
l\hbox{ if }i=0,
\end{cases}
\]
for ${\bf a}$ and the condition $\rho_1=0$ for $(\lambda,\rho )$. Let
$\overline{\bf a}$ and $(\overline\lambda,\overline\rho)$ be the configuration
and the corresponding rigged partition obtained from ${\bf a}$ and 
$(\lambda,\rho)$ by removing the weight $l$ particle to the far right.
Since the energy shift
when a weight $l$ particle passing the configuration $\overline{\bf a}$ is
given by $\sum_{i\geq2}A_{l,\lambda_i}$, and this is exactly the difference of
$\rho_1$ and $s_1$ in (\ref{RIGGING}) and (\ref{SR}), the above equivalence
follows.
\end{proof}

\section{Polynomial characters}
The purpose of this section is to derive fermionic character formulas for
the set of configurations with initial and boundary conditions.

We consider the $(k,3)$-configurations.
The initial conditions are specified by two integers $a$ and $b$
such that $0\leq a,b\leq k$: we set
\begin{equation}\label{AB}
C^{(k,l)}_{a,b}=\{{\bf a}\in C^{(k,l)}_{\rm pos};a_0=a,a_1=b\}.
\end{equation}
The problem is to determine the image of this set by the mapping $\iota$.
Note that (\ref{AB}) is empty unless $a+b\leq l$.

For a sequence of non-negative integers ${\bf r}=(r_1,\ldots,r_l)$
we define a subset $R^{(l)}({\bf r})$ of $R^{(l)}_{\rm pos}$ as follows.
\begin{equation}
R^{(l)}({\bf r})=\{(\lambda,\rho)\in R^{(l)};\rho_i\geq r_{\lambda_i}
\hbox{ for all }i\}.
\end{equation}
Let $J$ be a subset of $I=\{1,\ldots,l\}$. We define
\[
R^{(l)}({\bf r})_J=R^{(l)}({\bf r})\backslash R^{(l)}({\bf r}(J)),
\]
where
\[
{\bf r}(J)_i=
\begin{cases}
r_i+1&\hbox{if }i\in J;\\
r_i&\hbox{otherwise}.
\end{cases}
\]
In general, for a sequence of nonempty subsets
$J_m\subset I$ $(1\leq m\leq n)$, we set
\[
R^{(l)}({\bf r})_{J_1,\ldots,J_n}=R^{(l)}({\bf r})\backslash\left(
\bigcup_{m=1}^nR^{(l)}({\bf r}(J_m))\right).
\]
We have
\[
R^{(l)}({\bf r})=R^{(l)}({\bf r}(J))\bigsqcup R^{(l)}({\bf r})_J.
\]
Suppose that $J_1,J_2,J_3,\ldots,J_n\subset I$ are such that
$(J_1\cup J_2)\cap(\cup_{m=3}^nJ_m)=\emptyset$.
The following equalities are clear by the definition.
\begin{eqnarray}
&&R^{(l)}({\bf r})_{J'}=
R^{(l)}({\bf r})_{J_1,J'}\sqcup R^{(l)}({\bf r}(J_1))_{J'},\label{EQ}\\
&&R^{(l)}({\bf r})_{J_1\cup J_2,J'}=
R^{(l)}({\bf r})_{J_1,J'}\sqcup R^{(l)}({\bf r}(J_1))_{J_2,J'}
\hbox{ if }J_1\cap J_2=\emptyset,\label{EQA}\\
&&R^{(l)}({\bf r})_{J_1,J_2,J'}=
R^{(l)}({\bf r})_{J_1\cap J_2,J'}\bigsqcup
R^{(l)}({\bf r}(J_1\cap J_2))_{J_1\backslash J_2,J_2\backslash J_1,J'},
\label{EQB}\end{eqnarray}
where we used $J'$ to mean $J_3,\ldots,J_n$ for notational simplicity.
For example,
\[
R^{(l)}({\bf r})_{J'}=R^{(l)}({\bf r})_{J_3,\ldots,J_n}.
\]
If $J_1\subset J_2\subset I$, we have
\begin{equation}
R^{(l)}({\bf r})_{J_1,J_2}=R^{(l)}({\bf r})_{J_1}.
\end{equation}

\begin{thm}\label{INIT}
Let $[a,b]_l$ $(0\leq a,b\leq l;a+b\leq l)$
be the image of $C^{(k,l)}_{a,b}$ by the mapping $\iota$.
This is independent of $k$, and is given by
\begin{equation}\label{SETS}
[a,b]_l=
\begin{cases}
R^{(l)}({\bf r}_{a,b})_{[a,a+b],[a+b,l]}&\hbox{if }a\not=0,\\
R^{(l)}({\bf r}_{a,b})_{[b,l]}&\hbox{if }a=0\hbox{ and }b\not=0,\\
R^{(l)}({\bf r}_{a,b})&\hbox{if }a=0\hbox{ and }b=0,
\end{cases}
\end{equation}
where
\[
{\bf r}_{a,b}=(\underbrace{0,\ldots,0}_a,\underbrace{1,\ldots,b}_b,
\underbrace{b+2,\ldots,2l-2a-b}_{l-a-b}),
\]
and $[l_1,l_2]=\{l_1,l_1+1\ldots,l_2\}$ for
$1\leq l_1\leq l_2\leq l$.
\end{thm}

For a subset $R$ of $R^{(k)}$ we denote by $\chi(R)$ its character
\[
\chi(R)=\sum_{(\lambda,\rho)\in R}q^{d(\lambda,\rho)},
\]
where $d(\lambda,\rho)$ is given by (\ref{DEG}).

\begin{cor}
We have the following identities for the characters.
\begin{eqnarray}
\chi([a,b]_l)=
&&\chi(R^{(l)}(\underbrace{0,\ldots,0}_a,\underbrace{1,\ldots,b}_b,
\underbrace{b+2,\ldots,2l-2a-b}_{l-a-b}))\nonumber\\
&&-\chi(R^{(l)}(\underbrace{0,\ldots,0}_{a-1},1,\underbrace{2,\ldots,b+1}_b,
\underbrace{b+2,\ldots,2l-2a-b}_{l-a-b}))\nonumber\\
&&-\chi(R^{(l)}(\underbrace{0,\ldots,0}_a,\underbrace{1,\ldots,b-1}_{b-1},b+1,
\underbrace{b+3,\ldots,2l-2a-b+1}_{l-a-b}))\nonumber\\
&&+\chi(R^{(l)}(\underbrace{0,\ldots,0}_{a-1},1,
\underbrace{2,\ldots,b}_{b-1},b+1,
\underbrace{b+3,\ldots,2l-2a-b+1}_{l-a-b})),
\end{eqnarray}
where terms with $\underbrace{*}_{-1}$ is understood as $0$.
The second and the last term cancels each other if $a+b=l$ except for $a=l$.
\end{cor}

We prove this theorem by induction on $l$. In the following we abbreviate
$R^{(l)}({\bf r})$ to $({\bf r})$. For $l=1$, the statement of the
theorem is that
\[
[0,0]_1=(2),[0,1]_1=(1)\backslash(2),[1,0]_1=(0)\backslash(1).
\]
This is obvious because for ${\bf a}\in C^{(k,1)}$ the lowest position
of the (weight $1$) particles in ${\bf a}$ is equal to $\rho_m$ where
$\iota({\bf a})=((1^m),\rho)$. This is the base of the induction.

Note also that the theorem implies that the first two elements in the
configuration ${\bf a}=(a_0,a_1,\ldots)\in C^{(k,l)}$, i.e., $a_0$ and $a_1$,
are uniquely determined by the set of integers
$\rho^{(i)}_{m_i}$ $(1\leq i\leq l)$, where we use the notation
$\lambda=(\underbrace{l,\ldots,l}_{m_l},\ldots,\underbrace{1,\ldots,1}_{m_1})$
and $\rho=(\rho^{(l)}_1,\ldots,\rho^{(l)}_{m_l},\ldots,
\rho^{(1)}_1,\ldots,\rho^{(1)}_{m_1})$ for
${\bf a}\in C^{(k,l)}$ with $\iota({\bf a})=(\lambda,\rho)$, and we set
formally $\rho^{(i)}_0=\infty$ when $m_i=0$. This statement is also obvious
by the following reason. Since our construction of the bijections
proceed inductively on $l$, it is enough to show this statement for $i=l$.
By Proposition \ref{LPC} the values $\rho^{(l)}_j$ for $j<m_l$
are uniquely determined by $(a_2,a_3,\ldots)$. Conversely, the position of the
second lowest particles does not effect the values of $a_0$ and $a_1$.

By the definition the subsets $[a,b]_l$ are disjoint and the union
is equal to $R^{(l)}$. Therefore, the subsets in the right hand side
of (\ref{SETS}) must enjoy the same property. This statement will be directly
checked in the proof of the theorem. In the following we use the notation $A+B$
to mean the union of $A$ and $B$, and that $A$ and $B$ are disjoint.

The induction proceeds by the following recursion relation for the subsets
$[a,b]_l$. We define operations of constructing a subset of $R^{(l)}$
out of a subset $U$ of $R^{(l-1)}$. For $(\lambda,\rho)\in R^{(l)}$
we set $(\bar\lambda,\bar\rho)\in R^{(l-1)}$ by dropping the parts $\lambda_i$
(and $\rho_i$) such that $\lambda_i=l$. We set
\[
U*c=\{(\lambda,\rho)\in R^{(l)};
(\bar\lambda,\bar\rho)\in U,\rho^{(l)}_{m_l}\geq c\},
\quad U*{\underline c}=(U*c)\backslash(U*(c+1)).
\]
\begin{prop}
The subsets $[a,b]_l$ are determined by the following recursion relations.
\begin{equation}\label{REC}
[a,b]_l=
\begin{cases}
[a,b]_{l-1}*(2l-2a-b)+\sum_{c=0}^{b-1}
\bigl([a,c]_{l-1}*\underline{(2l-2a-b)}\bigr)
&\hbox{if }a+b<l;\\[5pt]
\sum_{c=0}^a\sum_{d=0}^{l-c-1}
\bigl([c,d]_{l-1}*\underline{(l-a)}\bigr)&\hbox{if }a+b=l.
\end{cases}
\end{equation}
\end{prop}
\begin{proof}
We use the notation $(\bar\lambda,\bar\rho)$ as above.
Set $\bar{\bf a}=\kappa(\bar\lambda,\bar\rho)$. Suppose that
\[
\bar{\bf a}=(c,d,\ldots).
\]
As we have remarked, it is enough to consider the configurations
where the number of weight $l$ particles is $1$. Consider the left moves
of a configuration ${\bf b}$ obtained from $\bar{\bf a}$ by adding a weight $l$
particle at a sufficiently large energy. We know that in finite,
say $t$, steps, the weight $l$ particle reaches the energy $0$.
Namely, the configuration $M_-^{t_0}{\bf b}$ is such that $(l,0,\ldots)$.
The values of the rigging corresponding to this particle is $0$ by the
definition. Let us consider how the configuration changes from 
$(c,d,\ldots)$ to $(l,0,\ldots)$. The change from $M_-^t{\bf b}$ to
$M_-^{t+1}{\bf b}$ is such that $+1$ at a column, say the $i$-th column,
and $-1$ at the $(i+1)$-th column. We have $i=0$ if and only if
$(M_-^t{\bf b})_0+(M_-^t{\bf b})_1=l$ since $(M_-^t{\bf b})_{-1}=0$.
Therefore, we have
\[
((M_-^t{\bf b})_0,(M_-^t{\bf b})_1)=
\begin{cases}
(c,d)&\hbox{if }t\leq t_0-2l+2c+d;\\
(c,t-t_0+2l-2c)&\hbox{if }t_0-2l+2c+d\leq t\leq t_0-l+c;\\
(l+t-t_0,t_0-t)&\hbox{if }t_0-l+c\leq t.
\end{cases}
\]
Therefore, a configuration ${\bf a}=(a,b,\ldots)$ appears
in this sequence if and only if
\[
\begin{cases}
a=c\hbox{ and }b\geq d&\hbox{for }a+b<l;\\
a\geq c&\hbox{for }a+b=l.
\end{cases}
\]
Counting the number of steps for $(a,b,\ldots)$
to reach $(l,0,\ldots)$, we obtain the value of the rigging.
If $a+b<l$, the change of the first two columns is such that
\[
(a,b)\rightarrow(a,b+1)\rightarrow\cdots\rightarrow(a,l-a)\rightarrow
(a+1,l-a-1)\rightarrow\cdots\rightarrow(l,0).
\]
Namely, $(a,b)$ reaches $(l,0)$, where the rigging is $0$, by
$2l-2a-b$ steps. This observation gives
(\ref{REC}) in the case $a+b<l$. The case $a+b=l$ is similar.
\end{proof}
{\it Proofs of Theorem \ref{INIT}.}
It is enough to show that the right hand side of (\ref{SETS}) satisfies
the recursion relation (\ref{REC}).

{\it Case $a+b<l$.}
Recall that we abbreviate $R^{(l)}({\bf r})$ to $({\bf r})$.
We first sum $[a,0]_{l-1}$ with $[a,1]_{l-1}$.
By using (\ref{EQB}), we decomposes $[a,1]_{l-1}$.
\begin{eqnarray}
&&(\underbrace{0,\ldots,0}_a,1,3,\ldots,2l-2a-3)_{[a,a+1],[a+1,l-1]}\nonumber\\
&&=(\underbrace{0,\ldots,0}_a,1,3,\ldots,2l-2a-3)_{[a+1,a+1]}
+(\underbrace{0,\ldots,0}_a,2,3,\ldots,2l-2a-3)_{[a,a],[a+2,l-1]}.\nonumber\\
\label{SUM1}
\end{eqnarray}
By using (\ref{EQ}), we sum the second term in the right hand side
with $[a,0]_{l-1}$.
\begin{eqnarray*}
&&(\underbrace{0,\ldots,0}_a,2,4,\ldots,2l-2a)_{[a,a]}+
(\underbrace{0,\ldots,0}_a,2,3,\ldots,2l-2a-3)_{[a,a],[a+2,l-1]}\\
&&=(\underbrace{0,\ldots,0}_a,2,3,\ldots,2l-2a-3)_{[a,a]}.
\end{eqnarray*}
By using (\ref{EQA}), we sum the first term in the right hand side of
(\ref{SUM1}) with this.
\begin{eqnarray*}
&&(\underbrace{0,\ldots,0}_a,1,3,\ldots,2l-2a-3)_{[a+1,a+1]}
+(\underbrace{0,\ldots,0}_a,2,3,\ldots,2l-2a-3)_{[a,a]}\\
&&=(\underbrace{0,\ldots,0}_a,1,3,5,\ldots,2l-2a-3)_{[a,a+1]}.
\end{eqnarray*}
Therefore, we obtain
\[
[a,0]_{l-1}+[a,1]_{l-1}=
(\underbrace{0,\ldots,0}_a,1,3,5,\ldots,2l-2a-3)_{[a,a+1]}.
\]
We repeat a similar summation until we obtain
\begin{eqnarray}
&&[a,0]_{l-1}*\underline{(2l-2a-b)}+\cdots+[a,b-1]_{l-1}*\underline{(2l-2a-b)}
\nonumber\\[3pt]&&
=(\underbrace{0,\ldots,0}_a,1,2,\ldots,b-1,b+1,b+3,
\ldots,2l-2a-b-1,2l-2a-b)_{[a,a+b-1],[l,l]}.\nonumber\\[-10pt]\label{SUM2}
\end{eqnarray}
Finally, we sum this result with
\begin{eqnarray}
&&[a,b]_{l-1}*(2l-2a-b)\nonumber\\
&&=(\underbrace{0,\ldots,0}_a,1,2,\ldots,b-1,b,b+2,
\ldots,2l-2a-b)_{[a+b,a+b]}\nonumber\\
&&\hskip-15pt+(\underbrace{0,\ldots,0}_a,1,2,\ldots,b-1,b+1,b+2,
\ldots,2l-2a-b)_{[a,a+b-1],[a+b+1,l-1]}.\label{SUM3}
\end{eqnarray}
If $a+b=l-1$, we must drop
the second term from the right hand side of this identity.
By using (\ref{EQA}), we sum (\ref{SUM2}) with this second term and obtain
\[
(\underbrace{0,\ldots,0}_a,1,2,\ldots,b-1,b+1,b+2,
\ldots,2l-2a-b-2,2l-2a-b)_{[a,a+b-1],[a+b+1,l]}.
\]
We sum this result with the first term in the right hand side of (\ref{SUM3})
and obtain $[a,b]_l$.

{\it Case $a+b=l$.}
Similarly, if $c>0$, we have
\[
\sum_{d=0}^{l-c-2}[c,d]_{l-1}
=(\underbrace{0,\ldots,0}_c,1,2,\ldots,l-2-c,l-c)_{[c,l-2]}.
\]
By using (\ref{EQA}), we sum this result with
\[
[c,l-1-c]_{l-1}=
(\underbrace{0,\ldots,0}_c,1,2,\ldots,l-2-c,l-1-c)_{[l-1,l-1]},
\]
and obtain
\[
(\underbrace{0,\ldots,0}_c,1,2,\ldots,l-2-c,l-1-c)_{[c,l-1]}.
\]
By using (\ref{EQ}), we first obtain
\[
\sum_{d=0}^{l-2}[0,d]_{l-1}=(1,2,\ldots,l-1),
\]
and then obtain
\[
\sum_{c=0}^a\sum_{d=0}^{l-c-1}[c,d]_{l-1}=
(\underbrace{0,\ldots,0}_a,1,2,\ldots,l-1-a).
\]
We obtain (\ref{REC}) for $a+b=l$ from this.
\qed

Now, we consider configurations zero at the boundary, i.e.,
above certain energy level. Set
\begin{equation}
C^{(k,l)}_{\rm pos}[N]=\{{\bf a}\in C^{(k,l)}_{\rm pos};
a_i=0\hbox{ for all }i>N\}.
\end{equation}
The following theorem describes the image of this finite set
in $R^{(l)}$ by the bijection $\iota$.
\begin{thm}\label{N}
A configuration ${\bf a}$ belongs to $C^{(k,l)}_{\rm pos}[N]$
if and only if the corresponding rigged partition
$(\lambda,\rho)=\iota({\bf a})$ satisfies
\begin{equation}\label{IN}
\rho_i\leq \lambda_iN-\sum_{j\not=i}A_{\lambda_i,\lambda_j}.
\end{equation}
\end{thm}

To prove this theorem, we prepare a few lemmas. Proofs are straightforward.
\begin{lem}\label{OUT}
Suppose that ${\bf a}\in C^{(k,l)}_{\rm pos}$ and
$(\lambda,\rho)=\iota({\bf a})$. Let ${\bf b}$ be the configuration
obtained from ${\bf a}$ by the parallel shift$:$ $b_i=a_{i-1}$.
Set $(\lambda',\rho')=\iota({\bf b})$. Then, we have
\[
\lambda'_i=\lambda_i,\quad\rho'_i=\rho_i+\lambda_i.
\]
\end{lem}
\begin{lem}\label{LM}
Let $1\leq l'<l\leq k$.
Suppose ${\bf a}\in C^{(k,l')}$ is such that $a_i=0$ for all $i<0$.
Then, we have $P_l{\bf a}_i=0$ for all $i<2$.
\end{lem}

{\it Proof of Theorem \ref{N}.}
Let $\lambda=(\lambda_1,\ldots,\lambda_n)$ and $\rho=(\rho_1,\ldots,\rho_n)$.
We assume that $\lambda_1=l$. Suppose that ${\bf a}\in C^{(k,l)}_m$, i.e.,
$\lambda_1=\cdots=\lambda_m=l>\lambda_{m+1}$, where $1\leq m\leq n$.

{\it Proof of ``only if'' part}.
We use an induction on $l$.
If $l=0$, there is nothing to prove.
We assume that ${\bf a}\in C^{(k,l)}_{\rm pos}[N]$.
First we show that $\rho_1\leq lN-\sum_{i=2}^nA_{l,\lambda_i}$.
For some $t$ the right move ${\bf a}'=M_+^t{\bf a}$ becomes
\[
a'_i=
\begin{cases}
0&\hbox{if }i\geq N\hbox{ or }i<0;\\
l&\hbox{if }i=N.
\end{cases}
\]
By the definition of the mapping $\iota$ we have
\[
\rho_1=lN-t-\sum_{i=2}^nA_{l,\lambda_i}.
\]
Therefore, we have (\ref{IN}) for $i=1$.
Since $\lambda_i\leq\lambda_1$ for all $2\leq i\leq m$, we have
(\ref{IN}) for all $2\leq i\leq m$.

Now, we will show (\ref{IN}) for $m+1\leq i\leq n$. Recall the procedure
of finding $\lambda$ and $\rho$. We bring all the weight $l$ particles
in ${\bf a}$ to a free position by the right move
$(M^{(m)}_+)^t\cdots(M^{(1)}_+)^t$ for a sufficiently large $t$.
The rest of the configuration ${\bf a}''$ is independent of $t$,
and it is supported in the finite interval $\{0,\ldots,N-2m\}$.
By the definition
\[
((\lambda_{m+1},\ldots,\lambda_n),(\rho_m+1,\ldots,\rho_n))
=\iota({\bf a}'').
\]
Now, let the weight $l$ particles in 
$(M^{(m)}_+)^t\cdots(M^{(1)}_+)^t{\bf a}$
pass the configuration ${\bf a}''$ from the right to the left
one by one. The configuration ${\bf a}''$ belongs to $C^{(k,l-1)}$.
By Lemma \ref{LM}, the resulting configuration is
supported in the interval $\{2m,\ldots,N\}$. Using Proposition \ref{SHIFT}
and Lemma \ref{OUT}, and also the induction hypothesis, we obtain
\[
\rho_i+mA_{l,\lambda_i}-2m\lambda_i\leq\lambda_i(N-2m)
-\sum_{m+1\leq j\leq n\atop j\not=i}A_{\lambda_i\lambda_i}.
\]
This is nothing but (\ref{IN}) for $m+1\leq i\leq n$.

{\it Proof of ``if'' part}.
We use induction on $l$. Assume that
\begin{equation}\label{ASS}
\rho_i\leq
\begin{cases}
lN-(m-1)A_{l,l}-{\displaystyle\sum_{j=m+1}^n}
A_{l,\lambda_j}&\hbox{ if }i\leq m;\\[15pt]
\lambda_iN-mA_{\lambda_i,l}
-{\displaystyle\sum_{m+1\leq j\leq n\atop j\not=i}}A_{\lambda_i,\lambda_j}
&\hbox{ if }i\geq m+1.
\end{cases}
\end{equation}
We move the weight $l$ particles to the far left.
Denote by ${\bf a}'''$ the rest of the configuration.
By Proposition \ref{SHIFT}, we see that the particle content of ${\bf a}'''$
is $(\lambda_{m+1},\ldots,\lambda_n)$, and the rigging is
$(\rho_{m+1}+mA_{\lambda_i,l},\ldots,\rho_n+mA_{\lambda_i,l})$.
The assumption (\ref{ASS}) implies
\[
\rho_i+mA_{\lambda_i,l}\leq\lambda_iN
-{\displaystyle\sum_{m+1\leq j\leq n\atop j\not=i}}A_{\lambda_i,\lambda_j}.
\]
Therefore, by the induction hypothesis, we have
\begin{equation}\label{a'''}
a'''_i=0\hbox{ for all }i>N.
\end{equation}
Next, starting from the original configuration
${\bf a}$, we move the weight $l$ particle at the highest position to the far
right, say to the energy, say $d$. Because of (\ref{a'''}),
in the process of reaching the level $d$, this particle must go through
the energy $lN$. In other words, before this moment in the up-going
process the whole configuration is supported in the region $i\leq N$.
After that, the further move breaks the support condition, and the weight
$l$ particle reaches the energy $d$. By the definition of the rigging,
in order to get back to ${\bf a}$, we must move this particle
to the left by
\[
d-\rho_1-(m-1)A_{l,l}-{\displaystyle\sum_{j=m+1}^n}A_{l,\lambda_j}
\]
steps. From (\ref{ASS}), we see that the number of steps
is greater than or equal to $d-lN$. This implies the original configuration
${\bf a}$ belongs to $C^{(k,l)}_{\rm pos}[N]$.
\qed

The polynomial identities (\ref{POL}) follow from Theorems \ref{INIT} and
\ref{N}. Set
\[
C^{(k,l)}_{a,b}[N]=C^{(k,l)}_{\rm pos}[N]\cap C^{(k,l)}_{a,b}.
\]
\begin{thm}
Suppose that $N\geq0$, $1\leq l\leq k$, $0\leq a,b$ and $a+b\leq l$.
We have the following identities.
\begin{equation}
\sum_{{\bf a}\in C^{(k,l)}_{a,b}[N]}q^{d({\bf a})}
=
\begin{cases}
\chi^{(k,l)}_{a,b}[N]-\chi^{(k,l)}_{a-1,b+2}[N]-\chi^{(k,l)}_{a,b-1}[N]
+\chi^{(k,l)}_{a-1,b+1}[N],&\hbox{if }b>0;\\
\chi^{(k,l)}_{a,0}[N]-\chi^{(k,l)}_{a-1,2}[N]&\hbox{if }b=0,
\end{cases}
\end{equation}
where $\chi^{(k,l)}_{a,b}[N]$ is given by
\begin{equation}\label{LAST}
\sum_{m_1,\ldots,m_k=0}^\infty
q^{Q({\bf m})+\sum_{i=1}^kr_im_i}\prod_{1\leq j\leq k\atop m_j\not=0}
\left[jN-\sum_{i=1}^kA_{j,i}m_i+A_{j,j}-r_j+m_j\atop m_j\right]
\end{equation}
with
\begin{eqnarray*}
Q({\bf m})&=&\frac12(A{\bf m},{\bf m})-\frac12\sum_{j=1}^kA_{j,j}m_j,
\quad{\bf m}=(m_1,\ldots,m_k),\\
(r_1,\ldots,r_k)&=&
(\underbrace{0,\ldots,0}_a,\underbrace{1,\ldots,b}_b,
\underbrace{b+2,\ldots,2k-2a-b}_{k-a-b}),
\end{eqnarray*}
and the summation is restricted to $m_{l+1}=\cdots=m_k=0$.
We understand $\chi^{(k,l)}_{a-1,b+2}[N]=\chi^{(k,l)}_{a-1,k-a+1}[N]$
if $a+b=k$, and $\chi^{(k,l)}_{a,b}[N]=0$ if $a=-1$.
\end{thm}

\section{Appendix}
The appendix contains Lemmas used in the proof of Propositions \ref{L'}
and \ref{SHIFT}.

Fix $1\leq l'<l\leq k$.
Recall Definition \ref{DEF} of the set of configurations $C^{(k,l')}$.
We consider an element ${\bf a}$ in $C^{(k,l')}$, and
${\bf a}'=P{\bf a}\in C^{(k,l)}$ (see Proposition \ref{L'}).

In the below we use frequently the equality
$a_{i}+a_{i+1}''=a_{i}''+a_{i+1}'$.

\begin{lem}\label{LA}
Suppose that the history for ${\bf a}\rightarrow{\bf a}'$ contains
\begin{eqnarray*}
S_2&:&\ldots,0,a_1,a_2,a''_3,a'_4,\ldots,\\
S_1&:&\ldots,0,a_1,a''_2,a'_3,a'_4,\ldots.
\end{eqnarray*}
Then, we have $a'_3+a'_4\leq l'$.
\end{lem}
\begin{proof}
We have $a_1+2a_2+2a''_3+a'_4\leq k+l$.
Since $a_2+a''_3=l$, we have $a_1+a'_4\leq k-l\leq l'$.
Since $a_1=a'_3$ we have $a'_3+a'_4\leq l'$.
\end{proof}

\begin{lem}\label{LB}
Suppose that the  history for ${\bf a}\rightarrow{\bf a}'$ contains
\[
S_i:\ldots,a_{i-1},a_i,a''_{i+1},a'_{i+2},a'_{i+3},a'_{i+4}\ldots.
\]
Then, we have $a'_{i+1}+2a'_{i+2}+2a'_{i+3}+a'_{i+4}<k+l'$.
\end{lem}
\begin{proof}
We have $a_{i+1}+2a''_{i+2}+2a'_{i+3}+a'_{i+4}\leq k+l$.
Since $a_{i+1}+a''_{i+2}=a''_{i+1}+a'_{i+2}$, $a_i+a''_{i+1}=l$ and $a'_{i+1}<a''_{i+1}$, we have
$-a_{i+1}+(l-a_i)+a'_{i+1}+2a'_{i+2}+2a'_{i+3}+a'_{i+4}<k+l$.
Therefore, we have
$a'_{i+1}+2a'_{i+2}+2a'_{i+3}+a'_{i+4}<k+a_i+a_{i+1}\leq k+l'$.
\end{proof}

\begin{lem}\label{LCD}
Suppose that the history for ${\bf a}\rightarrow{\bf a}'$ contains
\[
L_i:\ldots,a_{i-1},a_i,a''_{i+1},a'_{i+2},a'_{i+3},\ldots.
\]
Then, we have $a'_{i+2}+a'_{i+3}\leq a_{i-1}+a_i$.
If, in addition, the history contains
\[
L_{i+1}:\ldots,a_{i-1},a_i,a_{i+1},a''_{i+2},a'_{i+3},\ldots,
\]
then we have $a'_{i+2}+a'_{i+3}=a_{i-1}+a_i$.
\end{lem}
\begin{proof}
If $L_i$, we obtain 
$a_{i-1}+2a_i+2a''_{i+1}+a'_{i+2}=k+l$
and
$a_i+2a''_{i+1}+2a'_{i+2}+a'_{i+3}\leq k+l$.
Therefore, we obtain $a'_{i+2}+a'_{i+3}\leq a_{i-1}+a_i$.
If $L_i$ and $L_{i+1}$, we have
\begin{eqnarray*}
a_i+2a_{i+1}+2a''_{i+2}+a'_{i+3}&=&k+l,\\
2(a''_{i+1}+a'_{i+2})&=&2(a_{i+1}+a''_{i+2}),\\
k+l&=&a_{i-1}+2a_i+2a''_{i+1}+a'_{i+2}.
\end{eqnarray*}
Summing up, we obtain $a'_{i+2}+a'_{i+3}=a_{i-1}+a_i$.
\end{proof}

\begin{lem}\label{LC}
Suppose that the history for ${\bf a}\rightarrow{\bf a}'$ contains
\begin{eqnarray*}
S_3&:&\ldots,0,a_1,a_2,a_3,a''_4,a'_5,a'_6,\ldots,\\
L_2&:&\ldots,0,a_1,a_2,a''_3,a'_4,a'_5,a'_6,\ldots,\\
S_1&:&\ldots,0,a_1,a''_2,a'_3,a'_4,a'_5,a'_6,\ldots.
\end{eqnarray*}
Then, we have
$a'_3+a'_4\leq a_1+a_2$ and $a'_4+2a'_5+a'_6<k$.
\end{lem}
\begin{proof}
We have $a'_3=a_2+a''_3-a''_2=a_2+a''_3-(l-a_1)$
and $a'_4\leq l-a''_3$. Therefore, we have $a'_3+a'_4\leq a_1+a_2$.

We have
\begin{eqnarray*}
a'_4+2a'_5+a'_6+l&=&a_3+a''_4+a'_4+2a'_5+a'_6\\
&<&a_3+2a''_4+2a'_5+a'_6\\
&\leq&k+l.
\end{eqnarray*}
Therefore, we have $a'_4+2a'_5+a'_6<k$.
\end{proof}
\begin{lem}\label{LD}
Suppose that the history for ${\bf a}\rightarrow{\bf a}'$ contains
\begin{eqnarray*}
L_3&:&\ldots,0,a_1,a_2,a_3,a''_4,a'_5,a'_6,\ldots,\\
L_2&:&\ldots,0,a_1,a_2,a''_3,a'_4,a'_5,a'_6,\ldots,\\
S_1&:&\ldots,0,a_1,a''_2,a'_3,a'_4,a'_5,a'_6,\ldots.
\end{eqnarray*}
Then, we have $a'_3+a'_4\leq a_1+a_2$ and $a'_3+2a'_4+2a'_5+a'_6\leq k+l'$.
\end{lem}
\begin{proof}
By Lemma \ref{LCD} we obtain $a'_5+a'_6\leq a_2+a_3$.
Since the change from $L_2$ to $S_1$ is a multiple of $(-1,+1)$
at the columns indexed with $2$ and $3$, we have $a_2+a''_3=a''_2+a'_3$.
Also, $S_1$ implies $a_1+a''_2=l$. On the other hand, we have
$a''_3+a'_4\leq l$. Therefore,
we have $a'_3+a'_4\leq a_2+a''_3-(l-a_1)+l-a''_3=a_1+a_2$,
and also $k+l=a_1+2a_2+2a''_3+a'_4=a_1+a''_2+a'_3+a_2+a''_3+a'_4
=l+a'_3+a_2+a''_3+a'_4$. Therefore, we have $a'_3+a'_4+a'_5+a'_6\leq
k-a_2-a''_3+a_2+a_3$. Using $a''_3\geq a_3$ we have
$a'_3+a'_4+a'_5+a'_6\leq k$. By Lemma \ref{LCD} we obtain
$a'_4+a'_5\leq a_1+a_2\leq l'$, and
the second assertion of the lemma follows.
\end{proof}
\begin{lem}\label{LEFG}
Suppose that the history for ${\bf a}\rightarrow{\bf a}'$ contains
\begin{eqnarray*}
L_3&:&\ldots,0,a_1,a_2,a_3,a''_4,a'_5,a'_6,a'_7,\ldots,\\
L_2&:&\ldots,0,a_1,a_2,a''_3,a'_4,a'_5,a'_6,a'_7,\ldots,\\
S_1&:&\ldots,0,a_1,a''_2,a'_3,a'_4,a'_5,a'_6,a'_7,\ldots.
\end{eqnarray*}
Then, we have $a'_4+2a'_5+2a'_6+a'_7\leq k+l'$.
\end{lem}
\begin{proof}
There are three cases: in addition, the history contains $S_4$, $L_4$ or $L_5$.

Case $S_4$.
We have
\begin{eqnarray*}
2a'_6+a'_7&\leq&k-l+a_4,\\
a'_4+a'_5&\leq&l',\\
a'_5&=&l-a''_4.
\end{eqnarray*}
Summing up, we obtain $a'_4+2a'_5+2a'_6+a'_7\leq k+l'+a_4-a''_4\leq k+l'$.

Case $L_4$.
By Lemma \ref{LCD}, we obtain
$a'_4+2a'_5+2a'_6+a'_7\leq a_1+2a_2+2a_3+a_4\leq k+l'$.

Case $L_5$.
We have $a''_4=a_4$. Then, we have
$a'_4+2a'_5+2a'_6+a'_7=a_4+2a'_5+2a'_6+a'_7+a'_4-a_4
\leq k+l+a'_4-a_4$. Since $a_3+a_4=a''_3+a'_4$, we have
$a'_4-a_4=2a_3+a_4-2a''_3-a'_4$. Therefore, we obtain
$a'_4+2a'_5+2a'_6+a'_7\leq k+l+(a_1+2a_2+2a_3+a_4)-(a_1+2a_2+2a''_3+a'_4)
=a_1+2a_2+2a_3+a_4\leq k+l'$.
\end{proof}

\begin{lem}\label{LH}
Suppose that the history for ${\bf a}\rightarrow{\bf a}'$ contains
\begin{eqnarray*}
L_4&:&\ldots,0,a_1,a_2,a_3,a_4,a''_5,a'_6,a'_7,\ldots,\\
L_2&:&\ldots,0,a_1,a_2,a_3,a'_4,a'_5,a'_6,a'_7,\ldots,\\
S_1&:&\ldots,0,a_1,a''_2,a'_3,a'_4,a'_5,a'_6,a'_7,\ldots.
\end{eqnarray*}
Then, we have $(\ref{S3}$--$\ref{L6})$.
\end{lem}
\begin{proof}
First, note that $a'_3=a_2+a_3-a''_2=a_1+a_2+a_3-l$. Therefore, we have
$a'_3+a'_4=a_3+a'_4+a'_3-a_3\leq l+a_1+a_2-l=a_1+a_2\leq l'$.
By Lemma \ref{LCD}, we obtain $a'_4+a'_5\leq a_1+a_2$.

We have $a_3+2a_4+2a''_5+a'_6=k+l$.
This implies
\begin{equation}\label{FIN}
a_3+2a_4+a''_5\geq k,
\end{equation}
and also $a'_5+a'_6=k+l-(a_3+a_4+a'_4+a''_5)
=a_1+2a_2+2a_3+a_4+a'_4-(a_3+2a_4+a'_4+a''_5)
\leq k+l'-(a_3+2a_4+a''_5)$.
Using (\ref{FIN}), we obtain $a'_5+a'_6\leq l'$.

We have $a'_3+2a'_4+2a'_5+a'_6=a_3+2a'_4+2a'_5+a'_6+a'_3-a_3\leq
k+l+a_2-(l-a_1)=k+a_1+a_2\leq k+l'$.

Note that
$a'_5-a''_5=a_4-a'_4=a_1+2a_2+2a_3+a_4-(a_1+2a_2+2a_3+a'_4)\leq l'-l$.
Using this, we have
$a'_4+2a'_5+2a'_6+a'_7\leq a_4+2a''_5+2a'_6+a'_7+a'_5-a''_5\leq k+l'$,
and $a'_5+2a'_6+2a'_7+a'_8=a''_5+2a'_6+2a'_7+a'_8+a'_5-a''_5\leq k+l'$.
\end{proof}

\begin{lem}\label{LEM0}
If the history contains $S_2$, then we have $a_0+2a_1+2a_2+a_3<k+l'$.
\end{lem}
\begin{proof}
Note that $a_2+a_3<l'$ by (\ref{LOW}). Then, we have 
$a_0+2a_1+2a_2+a_3=a_0+2a_1+2a_2+a''_3+(a_2+a_3)-(a_2+a''_3)<k+l+l'-l=k+l'$.
\end{proof}
\begin{lem}\label{LEM1}
Suppose that the history contains
\begin{eqnarray*}
S_{i+1}&:&\ldots,a_{i-1},a_i,a_{i+1},a''_{i+2},a'_{i+3},\ldots,\\
S_i&:&\ldots,a_{i-1},a_i,a''_{i+1},a'_{i+2},a'_{i+3},\ldots.\\
\end{eqnarray*}
Then, we have $a_{i-1}+a_i\leq k-l$.
\end{lem}
\begin{proof}
Note that $a''_{i+1}=l-a_i$ and $a'_{i+2}=a_i$.
We have $a_{i-1}+2a_i+2(l-a_i)+a_i\leq k+l$, i.e., $a_{i-1}+a_i\leq k-l$.
\end{proof}
\begin{lem}\label{LEM2}
The history does not contain $S_2\leftarrow S_3$:
\begin{eqnarray*}
S_3&:&\ldots,a_1,a_2,a_3,a''_4,a'_5,\ldots,\\
S_2&:&\ldots,a_1,a_2,a''_3,a'_4,a'_5,\ldots.
\end{eqnarray*}
\end{lem}
\begin{proof}
By Lemma \ref{LEM0} we have $a_0+2a_1+2a_2+a_3<k+l'$.
By Lemma \ref{LEM1} we have $a_1+a_2\leq k-l<l'$. 
This is a contradiction to (\ref{HIGH}).
\end{proof}
\begin{lem}\label{LEM3}
The history does not contain $S_2\leftarrow L_3$:
\begin{eqnarray*}
L_3&:&\ldots,a_1,a_2,a_3,a''_4,a'_5,\ldots,\\
S_2&:&\ldots,a_1,a_2,a''_3,a'_4,a'_5,\ldots.
\end{eqnarray*}
\end{lem}
\begin{proof}
By Lemma \ref{LEM0} we have $a_0+2a_1+2a_2+a_3<k+l'$.
Note that $a_2+a_3<l'$ and $a''_4+a'_5\leq l$.
We have $a_1+a_2=a_1+2a_2+2a''_3+a'_4-(a_2+a''_3)-(a''_3+a'_4)
\leq k+l-l-(a_3+a''_4)=a_2+2a_3+2a''_4+a'_5-l-(a_3+a''_4)<l'+l-l=l'$.
This is a contradiction to (\ref{HIGH}).
\end{proof}
\begin{lem}\label{LEM4}
Suppose that $a_1+a_2=l'$, then, the history does not contain
$L_3\leftarrow L_5$:
\begin{eqnarray*}
L_5&:&\ldots,a_1,a_2,a_3,a_4,a_5,a''_6,a'_7\ldots,\\
L_3&:&\ldots,a_1,a_2,a_3,a_4,a'_5,a'_6,a'_7\ldots.
\end{eqnarray*}
\end{lem}
\begin{proof}
By (\ref{LOW}) we have $a_1+2a_2+2a_3+a_4=l'+a_2+2a_3+a_4<k+l'$.
Therefore, we have $a_2+2a_3+a_4<k$. Since $a_2+2a_3+2a_4+a'_5=k+l$,
we have $a_4+a'_5>l$. This is a contradiction.
\end{proof}
\begin{lem}\label{LEM5}
The history does not contain $L_2\leftarrow L_3\leftarrow L_5$:
\begin{eqnarray*}
L_5&:&\ldots,a_1,a_2,a_3,a_4,a_5,a''_6,a'_7\ldots,\\
L_3&:&\ldots,a_1,a_2,a_3,a_4,a'_5,a'_6,a'_7\ldots,\\
L_2&:&\ldots,a_1,a_2,a''_3,a'_4,a'_5,a'_6,a'_7\ldots.
\end{eqnarray*}
\end{lem}
\begin{proof}
By Lemma \ref{LEM4}, if $a_1+a_2=l'$, then the history does not contain
$L_3\leftarrow L_5$. Therefore, becauase of (\ref{HIGH}), we can assume that
$a_0+2a_1+2a_2+a_3=k+l'$. Then, we have
$a''_3-a_3=a_0+2a_1+2a_2+a''_3-(a_0+2a_1+2a_2+a_3)\leq l-l'$.
On the other hand, since $a_3+a_4=a''_3+a'_4$ we have
$a''_3-a_3=a_1+2a_2+2a''_3+a'_4-(a_1+2a_2+2a_3+a_4)>l-l'$.
This is a contradiction.
\end{proof}
\begin{lem}\label{LEMA}
If the  history contains $S_3$, then we have $a_1+a_2=l'$.
\end{lem}
\begin{proof}
The condition $S_3$ implies $a_3+a''_4=l$ and
$a_1+2a_2+2a_3+a''_4<k+l$. Therefore, we have $a_1+2a_2+a_3<k$,
and hence $a_0+2a_1+2a_2+a_3<k+l'$. By (\ref{HIGH}) we have $a_1+a_2=l'$.
\end{proof}
\begin{lem}\label{LEMC}
If the history contains the nodes $L_1$ and $L_4$,
\begin{eqnarray*}
L_4&:&\ldots,a_0,a_1,a_2,a_3,a_4,a''_5,a'_6,\ldots,\\
L_1&:&\ldots,a_0,a_1,a''_2,a'_3,a'_4,a'_5,a'_6,\ldots,
\end{eqnarray*}
then we have $L[4,{\bf a}']=L[1,{\bf a}]$.
\end{lem}
\begin{proof}
We have $L[1,{\bf a}]-L[4,{\bf a}']=
L[1,{\bf a}]-(a_0+2a_1+2a''_2+a'_3)
+a_3+2a_4+2a''_5+a'_6-L[4,{\bf a}']
=2(a_2+a_3+a_4+a''_5)-2(a''_2+a'_3+a'_4+a'_5)=0$.
The last equality follows from the observation that the moves between $L_4$
and $L_1$ are the move of $1$ inside the columns $2$ to $5$.
\end{proof}
\begin{lem}\label{S1MUST}
Suppose that the  history contains $L_2\leftarrow S_3$.
It continues as $S_1\leftarrow L_2$ or $L_1\leftarrow L_2$ or
$L_0\leftarrow L_2$. If $L_1\leftarrow L_2$ we have $a_1+a''_2=l$ at $L_1$.
Similarly, if $L_0\leftarrow L_2$ we have $a_1+a'_2=l$ at $L_0$.
In all cases, we have $a'_3+a'_4=a_1+a_2=l'$.
\end{lem}
\begin{proof}
By Lemma \ref{LEMA} we have $a_1+a_2=l'$.

Suppose that the history goes as
\begin{eqnarray*}
S_3&:&\ldots,a_0,a_1,a_2,a_3,a''_4,\ldots,\\
L_2&:&\ldots,a_0,a_1,a_2,a''_3,a'_4,\ldots,\\
L_1&:&\ldots,a_0,a_1,a''_2,a'_3,a'_4,\ldots.
\end{eqnarray*}
We have $a_1+2a_2+2a''_3+a'_4=k+l$. Using $a_1+a_2=l'$ and $a''_3+a'_4=l$,
we have $a_2+a''_3=a''_2+a'_3=k-l'$. Therefore, we have
$a_1+a''_2=a_0+2a_1+2a''_2+a'_3-(a_0+a_1)-(a''_2+a'_3)\geq k+l-l'-(k-l')=l$.
Therefore, we have $a_1+a''_2=l$.

Suppose that the history goes as
\begin{eqnarray*}
S_3&:&\ldots,a_{-1},a_0,a_1,a_2,a_3,a''_4,\ldots,\\
L_2&:&\ldots,a_{-1},a_0,a_1,a_2,a''_3,a'_4,\ldots,\\
L_0&:&\ldots,a_{-1},a_0,a_1,a'_2,a'_3,a'_4,\ldots,
\end{eqnarray*}
and that $a_1+a'_2<l$. Since $a_{-1}+2a_0+2a_1+a'_2=k+l$,
we have $a_{-1}+2a_0+a_1>k$. Since $a_{-1}+2a_0+2a_1+a_2\leq k+l'$,
we have $a_1+a_2<l'$. This is a contradiction.

Now, we show that $a'_3+a'_4=a_1+a_2$ in all cases.
Setting $a''_2=a'_2$ in the last case, we have $a''_2=l-a_1$ in all cases.
Then, we have $a'_3+a'_4=a_2+a''_3-a''_2+a_3+a''_4-a''_3=a_1+a_2$.
\end{proof}
\begin{lem}\label{EXCEPTIONAL}
Suppose that the  history contains
$L_2\leftarrow S_3\leftarrow S_4$,
\begin{eqnarray*}
S_4&:&\ldots,a_0,a_1,a_2,a_3,a_4,a''_5,a'_6,\ldots,\\
S_3&:&\ldots,a_0,a_1,a_2,a_3,a''_4,a'_5,a'_6,\ldots,\\
L_2&:&\ldots,a_0,a_1,a_2,a''_3,a'_4,a'_5,a'_6,\ldots.
\end{eqnarray*}
Then, we have $a_2+a_3\leq k-l$. If $a_2+a_3=k-l$, we have $a'_4+a'_5=l'$.
\end{lem}
\begin{proof}
By Lemma \ref{LEMA} we have $a_1+a_2=l'$. 
We have $a'_4=a_3+a''_4-a''_3=l-(k-a_1-2a_2)$ and
$a'_5=a_4+a''_5-a''_4=a_3$. Therefore, we obtain
$a'_4+a'_5=l-k+a_1+2a_2+a_3=l-k+l'+a_2+a_3$. Since $a'_4+a'_5\leq l'$ we have
$a_2+a_3\leq k-l$. Moreover, if $a_2+a_3=k-l$, we have $a'_4+a'_5=l'$.
\end{proof}

\begin{lem}\label{LEMF}
Suppose that the  history for ${\bf a}\rightarrow{\bf a}'$
contains the sequence
\begin{eqnarray*}
L_3&:&\ldots,a_{-1},a_0,a_1,a_2,a_3,a''_4,a'_5,\ldots,\\
L_2&:&\ldots,a_{-1},a_0,a_1,a_2,a''_3,a'_4,a'_5,\ldots,\\
L_0&:&\ldots,a_{-1},a_0,a_1,a'_2,a'_3,a'_4,a'_5,\ldots.
\end{eqnarray*}
If $a_1+a_2=l'$, then we have $a_1+a'_2=l$.
If $a_0+2a_1+2a_2+a_3=k+l'$, then we have $a_0+2a_1+2a'_2+a'_3=k+l$.
\end{lem}
\begin{proof}
Assume that $a_1+a_2=l'$.
From $a_{-1}+2a_0+2a_1+a_2\leq k+l'$ we have $a_{-1}+2a_0+a_1\leq k$.
Therefore, we obtain $a_1+a'_2=k+l-(a_{-1}+2a_0+a_1)\geq l$, and $a_1+a'_2=l$.

Assume that $a_0+2a_1+2a_2+a_3=k+l'$.
If $a_0+2a_1+2a'_2+a'_3<k+l$, we have
$k+l-(k+l')>a_0+2a_1+2a'_2+a'_3-(a_0+2a_1+2a_2+a_3)
=a'_2-a_2+a''_3-a_3\geq a'_2-a_2$. On the other hand, we have
$k+l-(k+l')\leq a_{-1}+2a_0+2a_1+a'_2-(a_{-1}+2a_0+2a_1+a_2)=a'_2-a_2$.
This is a contradiction.
\end{proof}
\begin{lem}\label{LEMG}
Suppose that the  history for ${\bf a}\rightarrow{\bf a}'$
contains the sequence
\begin{eqnarray*}
S_4&:&\ldots,a_0,a_1,a_2,a_3,a_4,a''_5,\ldots,\\
L_3&:&\ldots,a_0,a_1,a_2,a_3,a''_4,a'_5,\ldots,\\
L_2&:&\ldots,a_0,a_1,a_2,a''_3,a'_4,a'_5,\ldots.
\end{eqnarray*}
Then, we have $S[4,{\bf a}']=S[1,{\bf a}]=l'$.
\end{lem}
\begin{proof}
By Lemma \ref{LCD} we have $a'_4+a'_5=a_1+a_2$.
Suppose that $a_1+a_2<l'$. Then, we have $a_0+2a_1+2a_2+a_3=k+l'$.
There are three cases: (i) $S_1\leftarrow L_2$, (ii) $L_1\leftarrow L_2$
and (iii) $L_0\leftarrow L_2$. 

We lead to a contradiction in all cases.
We obtain the following equalities successively:
\begin{eqnarray*}
&&a''_5=l-a_4,\\
&&a''_4=k-a_2-2a_3,\quad a'_5=l-k+a_2+2a_3,\\
&&a''_3=l-a_1-a_2+a_3,\quad a'_4=k-l+a_1-2a_3.
\end{eqnarray*}

Case (i).
We have $a'_3=a_2+a''_3-a''_2=l-a_1+a_3-(l-a_1)=a_3$.
Since $a_0+2a_1+2a''_2+a'_3\leq k+l$, we have $a_0+a_3=a_0+a'_3\leq k-l$.
This implies $2(a_1+a_2)=k+l'-(a_0+a_3)\geq l+l'$, which is a contradiction.

Case (ii). We have $a''_2=k+l-a_0-2a_1-(a_2+a''_3)=k-a_0-a_1-a_3$.
This implies $0\leq a''_2-a_2=k-a_0-a_1-a_2-a_3=k-(k+l')+a_1+a_2=a_1+a_2-l'$,
which is a contradiction.

Case (iii). We formally set $a''_2=a'_2$. Then, by Lemma \ref{LEMF}, we have 
$a_0+2a_1+2a''_2+a'_3=k+l$. We can follow the proof for Case (ii).
\end{proof}

\begin{lem}\label{LEMG'}
Suppose that the  history for ${\bf a}\rightarrow{\bf a}'$
contains the sequence
\begin{eqnarray*}
L_4&:&\ldots,a_0,a_1,a_2,a_3,a_4,a''_5,\ldots,\\
L_3&:&\ldots,a_0,a_1,a_2,a_3,a''_4,a'_5,\ldots,\\
L_2&:&\ldots,a_0,a_1,a_2,a''_3,a'_4,a'_5,\ldots.
\end{eqnarray*}
Then, we have $S[4,{\bf a}']=S[1,{\bf a}]=l'$ or $L[4,{\bf a}']=k+l'$.
\end{lem}
\begin{proof}
The proof goes exactly the same as Lemma \ref{LEMG} for
the first paragraph. Then, we continue as follows.

Case (i). We have $l-l'\geq a_0+2a_1+2a''_2+a'_3-(a_0+2a_1+2a_2+a_3)
=a''_2-a_2+a''_3-a_3\geq a_1+a''_2-(a_1+a_2)>l-l'$. This is a contradiction.

Case (ii). By Lemma \ref{LEMC}, we have $L[4,{\bf a}']=L[1,{\bf a}]=k+l'$.
\end{proof}
\begin{lem}\label{LEMI}
Suppose that the  history for ${\bf a}\rightarrow{\bf a}'$
contains the sequence $L_2\leftarrow L_4$. Then, we have
$L[4,{\bf a}']=k+l'$.
\end{lem}
\begin{proof}
We set $a''_5=a_5$ if the history contains $L_4\leftarrow L_6$.
There are three cases: (i) $S_1\leftarrow L_2$, (ii) $L_1\leftarrow L_2$
and (iii) $L_0\leftarrow L_2$. Setting $a''_2=a'_2$ in (iii), we have the
sequence
\begin{eqnarray*}
L_4&:&\ldots,a_{-1},a_0,a_1,a_2,a_3,a_4,a''_5,a'_6,\ldots,\\
L_2&:&\ldots,a_{-1},a_0,a_1,a_2,a_3,a'_4,a'_5,a'_6,\ldots,\\
S_1\hbox{ or }L_1\hbox{ or }L_0
&:&\ldots,a_{-1},a_0,a_1,a''_2,a'_3,a'_4,a'_5,a'_6,\ldots.
\end{eqnarray*}

Case (i). We have $a'_3=a_2+a_3-(l-a_1)$. Therefore, we have
$L[4,{\bf a}']=a_1+a_2+a_3-l+2a_4+2a''_5+a'_6=a_1+a_2+k$.
If $a_1+a_2=l'$ we have $L[4,{\bf a}']=k+l'$.

Suppose that $a_1+a_2<l'$. We have $L[1,{\bf a}]=k+l'$, and, therefore,
$l-l'\geq a_0+2a_1+2a''_2+a'_3-(a_0+2a_1+2a_2+a_3)=a''_2-a_2+a''_3-a_3
\geq a''_2-a_2=a_1+a''_2-(a_1+a_2)>l-l'$. This is a contradiction.

Case (ii). By Lemma \ref{LEMC} we have $L[4,{\bf a}']=L[1,{\bf a}]$.
Therefore, if $L[1,{\bf a}]=k+l'$, we have $L[4,{\bf a}']=k+l'$.
If $L[1,{\bf a}]<k+l'$ we have $a_1+a_2=l'$. Therefore we have
$a''_2-a_2=a_1+a''_2-(a_1+a_2)\leq l-l'$. Since $a''_2+a'_3=a_2+a_3$, we have
$a''_2-a_2=a_0+2a_1+2a''_2+a'_3-(a_0+2a_1+2a_2+a_3)>l-l'$.
This is a contradiction.

Case (iii). We have $a''_2=k+l-(a_{-1}+2a_0+2a_1)$, and, therefore,
$a'_3=L[0,{\bf a}]+a_3-(k+l)$. This implies
$L[4,{\bf a}']=L[0,{\bf a}]+a_3-(k+l)+2a_4+2a''_5+a'_6=
L[0,{\bf a}]$. We will show that $L[0,{\bf a}]=k+l'$.

We have
\begin{equation}\label{LINE}
a''_2-a_2=a_{-1}+2a_0+2a_1+a''_2-(a_{-1}+2a_0+2a_1+a_2)\geq l-l'.
\end{equation}

If $L[1,{\bf a}]=k+l'$, we have
$a''_2-a_2=a_0+2a_1+2a''_2+a'_3-(a_0+2a_1+2a_2+a_3)\leq l-l'$.
Otherwise, we have $S[1, {\bf a}]=l'$ and
$a''_2-a_2=a_1+a''_2-(a_1+a_2)\leq l-l'$.
Therefore, in both cases, we have the equality at the end of (\ref{LINE}),
and, in particular, we have $L[0,{\bf a}]=k+l'$.
\end{proof}

\begin{lem}\label{LEMK}
Suppose that the  history for ${\bf a}\rightarrow{\bf a}'$
contains the sequence
\begin{eqnarray*}
S_4&:&\ldots,a_0,a_1,a_2,a_3,a_4,a''_5,\ldots,\\
L_3&:&\ldots,a_0,a_1,a_2,a_3,a''_4,a'_5,\ldots,\\
L_1&:&\ldots,a_0,a_1,a_2,a'_3,a'_4,a'_5,\ldots.
\end{eqnarray*}
Then, we have $S[4,{\bf a}']=L[1,{\bf a}]-k=l'$.
\end{lem}
\begin{proof}
We have the following equalities.
\begin{eqnarray*}
&&a'_3=k+l-(a_0+2a_1+2a_2),\\
&&a''_4=k-a_2-2a_3,\\
&&a'_4=a_0+2a_1+a_2-a_3-l,\\
&&a'_5=l-k+a_2+2a_3.
\end{eqnarray*}
Therefore, we have $S[4,{\bf a}']=L[1,{\bf a}]-k$.

If $L[1,{\bf a}]=k+l'$, the proof is over. Otherwise, we have $a_1+a_2=l'$.
Then, we have $k+l\geq a_1+2a_2+2a'_3+a'_4=2k+l-(a_0+a_1+a_2+a_3)$.
therefore, we have $L[1,{\bf a}]\geq k+l'$. This is a contradiction.
\end{proof}
\begin{lem}\label{LEML}
Suppose that the  history for ${\bf a}\rightarrow{\bf a}'$
contains the sequence
\begin{eqnarray*}
L_4&:&\ldots,a_0,a_1,a_2,a_3,a_4,a''_5,\ldots,\\
L_3&:&\ldots,a_0,a_1,a_2,a_3,a''_4,a'_5,\ldots,\\
L_1&:&\ldots,a_0,a_1,a_2,a'_3,a'_4,a'_5,\ldots.
\end{eqnarray*}
Then, we have $L[4,{\bf a}']=L[1,{\bf a}]=k+l'$ or
$S[4,{\bf a}']=S[1,{\bf a}]=l'$.
\end{lem}
\begin{proof}
By Lemma \ref{LEMC} we have $L[4,{\bf a}']=L[1,{\bf a}]$. If
$L[1,{\bf a}]=k+l'$, the proof is over. Otherwise, we have $a_1+a_2=l'$.
Then, we have $k+l\geq(a_1+a_2)+a_2+2a'_3+a'_4=l'+k+l-(a'_4+a'_5)$.
Therefore, we have $a'_4+a'_5\geq l'$. This implies $S[4,{\bf a}']=l'$.
\end{proof}
\begin{lem}\label{LEMM}
Suppose that the  history for ${\bf a}\rightarrow{\bf a}'$
contains the sequence
\begin{eqnarray*}
L_5&:&\ldots,a_0,a_1,a_2,a_3,a_4,a_5,a''_6,a'_7,\ldots,\\
L_3&:&\ldots,a_0,a_1,a_2,a_3,a_4,a'_5,a'_6,a'_7,\ldots,\\
L_1&:&\ldots,a_0,a_1,a_2,a'_3,a'_4,a'_5,a'_6,a'_7,\ldots.
\end{eqnarray*}
Then, we have $L[5,{\bf a}']=k+l'$.
\end{lem}
\begin{proof}
By Lemma \ref{LEM4} we have $L[1,{\bf a}]=k+l'$.
Therefore, we have $a'_3-a_3=a_0+2a_1+2a_2+a'_3-L[1,{\bf a}]=l-l'$.
Hence we have $a'_4-a_4=l'-l$. This implies
$L[5,{\bf a}']=a_4+2a_5+2a''_6+a'_7+(a'_4-a_4)=k+l+l'-l=k+l'$.
\end{proof}

\noindent
{\it Acknowledgments.}\quad 
B.F. is partially supported by the grants,
CRDF RP1-2254, INTAS 00-55, RFBR 00-15-96579.
M.J. is partially supported by
the Grant-in-Aid for Scientific Research (B2)
no.14340040, and 
T.M. is partially supported by
the Grant-in-Aid for Scientific Research (A)
no. 13304010, Japan Society for the Promotion of Science.
E.M. is partially supported by NSF grant DMS-0140460.


\end{document}